\documentclass[authoryear,12pt,3p]{elsarticle}
\usepackage[utf8]{inputenc}
\usepackage{amsmath}
\usepackage{amssymb}
\usepackage{multirow}
\usepackage{hyperref}
\usepackage{tikz}
\usepackage{tikz-network}
\usepackage{marvosym}
\usepackage[title]{appendix}
\usepackage{float}
\usepackage[ruled, lined]{algorithm2e}
\usepackage{setspace}
\usepackage[pagewise]{lineno}
\usepackage{enumitem}
\usepackage{xcolor}
\usepackage{caption}
\usepackage{subcaption}
\usetikzlibrary{shapes.geometric, arrows, decorations.pathreplacing,positioning, arrows.meta}

\newcommand{\ImageWidth}{11cm}

\definecolor{busabold}{HTML}{0606c9}
\definecolor{busafaded}{HTML}{9fb7d4}
\definecolor{busbfaded}{HTML}{e89309}
\definecolor{busbfaded}{HTML}{fdcd8a}  

\definecolor{fcblue}{rgb}{0.0, 0.18, 0.65}
\definecolor{tlblue}{rgb}{0.05, 0.1, 0.7}
\definecolor{fcgreen}{rgb}{0.0, 0.65, 0.18}
\definecolor{tlgreen}{HTML}{229936}
\definecolor{fcred}{rgb}{0.65, 0.05, 0.05}
\definecolor{pastelgreen}{HTML}{6ab55e}

\tikzstyle{start} = [rectangle, rounded corners, 
minimum width=3cm, 
minimum height=1cm,
text centered, 
draw=black,
text=white,
fill=fcgreen]

\tikzstyle{stop} = [rectangle, rounded corners, 
minimum width=3cm, 
minimum height=1cm,
text centered, 
draw=black,
text=white,
fill=fcred]

\tikzstyle{io} = [trapezium, 
trapezium stretches=true,
trapezium left angle=70, 
trapezium right angle=110, 
minimum width=3cm,
minimum height=1cm, text centered,
text=black,
draw=black, fill=gray!70]

\tikzstyle{process} = [rectangle, 
minimum width=3cm, 
minimum height=1cm, 
text centered, 
text width=3cm, 
draw=black,
text=white,
fill=fcblue]

\tikzstyle{decision} = [diamond, 
minimum width=3cm, 
minimum height=1cm,
aspect=2,
text centered, 
draw=black, 
fill=gray!70]
\tikzstyle{arrow} = [thick,->,>=stealth]

\allowdisplaybreaks
\newcommand{\setC}{\mathcal{C}}
\newcommand{\setT}{\mathcal{T}}
\newcommand{\setV}{\mathcal{V}}
\newcommand{\setA}{\mathcal{A}}
\newcommand{\setM}{\mathcal{M}}

\begin{document}
\begin{frontmatter}
\title{Minimum-Delay Opportunity Charging Scheduling for Electric Buses}
\author[1]{Dan McCabe} 
\ead{dmccabe@uw.edu}
\affiliation[1]{
    organization={Department of Civil and Environmental Engineering, University of Washington},
    city={Seattle},
    state={WA},
    country={USA}}

\author[1]{Xuegang (Jeff) Ban}
\ead{banx@uw.edu}

\author[2]{Bal\'azs Kulcs\'ar \corref{bk}}
\ead{kulcsar@chalmers.se}
\cortext[bk]{Corresponding author}

\affiliation[2]{
    organization={Department of Electrical Engineering, Chalmers University of Technology},
    city={Gothenburg},
    country={Sweden}}

\begin{abstract}
    \textcolor{black}{Transit agencies that operate battery-electric buses must carefully manage fast-charging infrastructure to extend daily bus range without degrading on-time performance. To support this need, we propose a mixed-integer linear programming model to schedule opportunity charging that minimizes the amount of departure delay in all trips served by electric buses. Our novel approach directly tracks queuing at chargers in order to set and propagate departure delays. Allowing but minimizing delays makes it possible to optimize performance when delays due to traffic conditions and charging needs are inevitable, in contrast with existing methods that require charging to occur during scheduled layover time. To solve the model, we develop two algorithms based on decomposition. The first is an exact solution method based on Combinatorial Benders (CB) decomposition, which avoids directly enumerating the model's logic-based ``big M'' constraints and their inevitable computational challenges. The second, inspired by the CB approach but more efficient, is a polynomial-time heuristic based on linear programming that we call 3S. Computational experiments on both a simple notional transit network and the real bus system of King County, Washington, USA demonstrate the performance of both methods. The 3S method appears particularly promising for creating good charging schedules quickly at real-world scale.}
\end{abstract}



\end{frontmatter}

\section{Introduction}
\label{sec:introduction}
Battery-electric buses (BEBs) make up a significant and growing share of the global transit vehicle fleet. Over 60,000 such vehicles were sold worldwide in 2022, about 5\% of global bus sales \citep{IEAOutlook2023}. While China has been the world leader in electric bus adoption for many years, BEBs are beginning to see greater adoption worldwide. In the United States, the Bipartisan Infrastructure Law of 2021 allocated over \$5 billion to help agencies purchase low-emissions transit vehicles and charging infrastructure \citep{FTAInfraBill}. Worldwide, BloombergNEF projects that 50\% of buses will be battery-powered by 2032, a milestone passenger cars are not expected to reach for a further ten years \citep{BNEF2023}.

As this rapid transition commences, agencies who replace some or all of their conventionally fueled buses with BEBs face several planning challenges across various time scales. Some key long-term decisions are bus fleet composition over time in terms of fuel type and/or battery size \citep{Rogge2018, Pelletier2019} as well as charging infrastructure design such as the location, number, and power output of chargers \citep{mccabe_optimal_2023, An2020}. After these decisions have been made, transit agencies must determine how to efficiently operate BEBs each day while providing good passenger service and maintaining charged batteries on all vehicles. The preferred approach for many agencies is to rely on low-power overnight charging at bus depots in order to mimic traditional bus operating patterns in which refueling is not a significant concern. However, in many cases agencies have to use high-power chargers during the day to effectively extend bus ranges without utilizing a larger battery \citep{LATransitionPlan,MBTAWinter,MetroTransitionPlan}. This approach is commonly referred to as \textit{opportunity charging} and sometimes \textit{on-route} or \textit{layover charging}. Transit operators who rely on opportunity charging need to schedule charging for all buses each day to make efficient use of their limited charging resources.

Charging scheduling has often been approached as a joint problem with vehicle scheduling, in which each bus that operates on a given day is assigned to a sequence of trips typically referred to as a \textit{block} \citep{Perumal2022}. Although combining these phases can help develop vehicle schedules well suited to BEBs, it is important to note that bus schedules are usually determined only once every few months, but agencies need to optimize charging schedules repeatedly at a variety of time scales. Models focused primarily on determining vehicle schedules (such as \citet{LiuCeder2020}) or charging locations (such as \citet{mccabe_optimal_2023}) generally design charging schedules based on average or expected conditions over their longer-range planning horizons, but within that time frame, operating conditions can easily change enough to require charging schedule adjustments. For example, weather and ridership patterns significantly impact energy consumption, and predictions of these variables improve as each service day approaches. For optimal daily operations, transit agencies should adjust their planned charging schedules frequently to account for these variations. Agencies also may wish to revise charging schedules within a single day to respond to changes such as vehicle breakdowns and traffic congestion---for example, it may be best to postpone a scheduled charge when a bus is running behind schedule, allowing it to resume service sooner and minimize delays. A revised charging schedule can be used both to improve operations and to provide more accurate estimates of future trip delays to passengers.

Researchers have developed various approaches to charge scheduling over the previous several years, but this work describes a novel modeling approach to the problem. First, our model prioritizes the passenger experience, minimizing departure delays across all trips in a day. We do this by explicitly modeling the queuing process at charging stations to track exactly when each bus is ready to start all of its trips. Our model propagates any delays across trips served by the same bus, ensuring their cascading impacts are captured accurately. Second, our focus on accurate delay tracking makes the model applicable to nearly any operating conditions. Rather than constraining all trips to leave on time like most other works in the literature, our model handles scenarios where zero delay is not achievable. It produces actionable results under all conditions and extends naturally to robust or stochastic applications where energy consumption and traffic conditions make on-time departures impossible.


\textcolor{black}{
Our precise modeling approach produces some computational challenges. In particular, the mixed-integer linear programming model includes many binary variables and ``big-M'' constraints. To mitigate this issue, we develop both an exact solution method and a polynomial-time heuristic algorithm. Both approaches decompose the complete problem into a series of easier-to-solve problems with easily interpretable structure. The exact method uses Combinatorial Benders (CB) decomposition to split the problem into a master problem that contains all binary variables and a subproblem with all continuous variables, which reduces the burden of the big-M constraints. This restructuring is natural for our problem: the binary variables reflect high-level decisions about the ordering of charging throughout the day, whereas the continuous variables track the exact timing of events, including charging times and trip departures. The Combinatorial Benders approach exploits this structure in an iterative framework.
}

\textcolor{black}{Our Select—Sequence—Schedule (3S) heuristic algorithm was inspired by the CB approach and is based on a similarly explainable decomposition. As in Combinatorial Benders, we make binary decisions first and continuous decisions second. To expedite this process, 3S initially relaxes the queue tracking constraints that link different buses to each other, so that we can \textit{select} when charging occurs with a separate linear program for each bus. We then establish the \textit{sequence} in which buses visit each charger with a simple sorting operation and finally \textit{schedule} the exact charger plugin times, charging durations, and delays with the same approach as in the CB subproblem. The 3S method has polynomial-time complexity in the worst case, so it can be used to solve real-world instances quickly. Our experiments show that despite its lack of a performance guarantee, 3S reliably identifies good feasible solutions (and often optimal solutions) much faster than exact approaches do.}

In summary, the main contributions of this work are as follows:

\begin{itemize}
    \item A mixed-integer linear programming model for recharging scheduling that exactly tracks queuing behavior at chargers, propagates delays across trips completed by the same bus, and minimizes the total delay.
    \item An exact solution method based on Combinatorial Benders decomposition to solve this computationally demanding model.
    \item A polynomial-time heuristic method, motivated by the exact solution method, which generates good solutions quickly for complex real-world transit networks.
\end{itemize}

The remainder of this paper is organized as follows. Section \ref{sec:litreview} provides a review of relevant literature. Section \ref{sec:methodology} describes our methodology and mixed-integer programming formulation. Section \ref{sec:benders} describes our exact solution method based on Combinatorial Benders decomposition, while Section \ref{sec:heuristic} documents the related heuristic. Section \ref{sec:casestudy} applies our methods to both a simple notional transit network and the real Seattle-area bus network. Section \ref{sec:conclusion} summarizes findings and concludes the paper.

\section{Literature Review}
\subsection{Recharging Scheduling for BEBs}
\label{sec:litreview}
Recharging scheduling for BEBs connects to a larger transit planning process that designs routes, timetables, and vehicle and driver schedules. In traditional transit planning with diesel buses, these problems are addressed sequentially rather than simultaneously, as each step is already complex to optimize \citep{Ceder2007}. Replacing diesel buses with BEBs further complicates the traditional transit planning process by introducing new challenges related to capital investments in buses and chargers, location of chargers, possible schedule revisions, and charging scheduling \citep{Perumal2022}. Because these various stages of the BEB planning process have significant overlap, researchers have studied recharging scheduling on its own as well as its integration with other decisions.

We first summarize the literature on charging scheduling integrated with other decisions, typically related to either charging infrastructure or vehicle schedules. When combining charge scheduling with other decisions, it is common to make significant assumptions to limit computational complexity, and these assumptions limit the models' applicability to day-to-day operations. For example, \citet{esmaeilnejad_optimal_2023} proposed a stochastic optimization model to determine both charger locations and recharging schedules along a single bus route, where charging may take place at intermediate stops but is penalized based on passenger waiting costs. However, this model assumes buses will always charge to 100\% state of charge between a pair of consecutive trips, which is only appropriate for buses with small battery capacities. \citet{mccabe_optimal_2023} developed a deterministic model for simultaneously locating chargers and scheduling charger usage. Charger capacity is accounted for with predetermined ``conflict sets'' based on posted timetables, which could result in unexpected queuing at chargers if buses are delayed. \citet{gairola_optimization_2023} proposed a comprehensive robust model to minimize the costs of converting to BEBs by optimizing battery sizes, charging infrastructure, and charging schedules. Their model allowed recharging at terminals during scheduled layover time and handled capacity by discretizing time into unit-length intervals and ensuring the number of buses charging at a terminal did not exceed the number of chargers there, similar to \citet{mccabe_optimal_2023}. However, neither of these approaches extends well to daily operations where unusually high energy demand or exogenous delays might make it impossible for buses to stay on schedule and charge only within scheduled layover time.

Charging scheduling is also commonly incorporated into vehicle block design, which for diesel buses is classically addressed with the vehicle scheduling problem (VSP). A comprehensive review of BEB VSP publications can be found in \citet{Perumal2022}, which references 23 such works. Among these, only six considered partial recharging of batteries and only seven accommodated multiple vehicle types, while none include both of these elements simultaneously. Both of these attributes are essential for daily recharging scheduling---many transit agencies operate both standard-length and articulated buses, for example, and partial recharging is a natural strategy for large batteries---but present major computational challenges to incorporate in a model that also designs blocks. 

Recognizing the need for a precise charge scheduling model focused on daily operations, some publications treat infrastructure and vehicle schedules as fixed and focus purely on charging scheduling. \citet{abdelwahed_evaluating_2020} developed two related formulations for charge scheduling at terminals and found that an event-based discretization scheme was preferable to using uniform time intervals. Both models sought to minimize the total cost of electricity based on time-of-use (TOU) prices and constrained charging to take place during the scheduled layover time. \citet{liu_robust_2022} developed a robust model for scheduling daily recharging, both at the depot and at terminals. Their model allows flexible recharging durations and power levels, respects charging station capacity, and minimizes the total cost of energy based on TOU pricing. They consider uncertainty in buses' energy consumption, but assume that buses always complete all trips as scheduled and require charging to occur during scheduled layover time. Using a different modeling paradigm, \citet{lacombe_integrated_2023} studied an optimal control problem intended to minimize both energy costs and schedule deviations, where the bus schedule is encoded as a route-specific headway rather than a timetable. They developed a decomposition strategy based on Lagrangian relaxation and local heuristics to apply the method in practice. This approach is much more flexible than restricting charging to scheduled layover time, but is only appropriate for systems where maintaining a target headway is more important than matching an advertised schedule. Finally, we note that there are various works such as \citet{manzolli2022electric} and \citet{brinkel_comparative_2023} that focus exclusively on optimizing depot charging. This is a fundamentally different challenge from our setting that is focused on opportunity charging at terminals (possibly including, but not limited to, depots), so we do not review those in detail here.

\subsection{Combinatorial Benders Decomposition}
Most approaches to optimal charge scheduling result in a mixed-integer programming problem, making them difficult to solve. We handle this challenge using a tailored version of Benders decomposition, a classical algorithm for mixed-integer programming that was first proposed over 60 years ago \citep{Benders1962}. The central idea of this approach is to separate the problem into a master problem (MP) containing integer variables and a subproblem (SP) containing continuous variables. In each iteration, the MP is solved to obtain candidate values of all integer variables. These fixed values are an input to the SP, which verifies the feasibility and potential optimality of the candidate solution. If the SP determines that the current candidate solution cannot be an optimal solution, one or more Benders cuts are generated and added to the MP to exclude it (along with, ideally, many other solutions), then the MP is solved again. The procedure repeats until an optimal solution is found \citep{rahmaniani2017benders}. 

Benders decomposition is often used in stochastic programming applications, where after fixing a limited number of first-stage variables, the second-stage problem decomposes into several independent problems that can be solved quickly. However, in the last two decades researchers have identified additional classes of problems where the approach can outperform standard branch-and-bound or branch-and-cut algorithms. \citet{Hooker2003} introduced the idea of logic-based Benders decomposition, wherein cuts are generated not based on the linear programming dual of the subproblem but a so-called ``inference dual'' that generalizes LP duality. Their approach has been applied successfully to a variety of problems but appears especially effective in cases where specialized constraint programming methods can be applied to the subproblem, including some types of scheduling problems \citep{hooker2007planning}.

\citet{Codato2006} soon after developed the idea of Combinatorial Benders (CB) cuts for a specific class of mixed-integer linear programs with logical constraints. Their approach largely follows \citet{Hooker2003}, but they derived a tailored method for finding cuts for their particular problem class and demonstrated its performance benefits on some example problems. A major selling point of this method is that it eliminates the computational problems caused by big-M constraints; once variable values have been set by the MP, the corresponding logical constraints are either included or excluded from the SP. As such, the CB approach can avoid the usual problem of big-M values giving a poor linear programming relaxation and resulting bounds.

It should be noted that a direct application of the Benders (standard or combinatorial) algorithm can have poor performance for a variety of reasons. These issues and strategies to mitigate them are reviewed thoroughly in \citet{rahmaniani2017benders}. Among the keys to a successful Benders implementation are initializing the MP with a set of strong cuts to aid in finding feasible solutions, using heuristics to generate good solutions and their correspondingly strong Benders cuts, and embedding Benders cuts within a branch-and-cut algorithm to limit the amount of redundant computations. We found these strategies were critical to make the CB approach competitive with an off-the-shelf solver.

\subsection{Summary}
The worldwide growth of BEBs has produced a significant body of literature on charging scheduling. However, most of these approaches are more suited to long-range planning than day-to-day or real-time scheduling. When charging scheduling is embedded within a model primarily focused on designing a charging infrastructure network or vehicle blocks, it tends to be simplified and restrictive. Even models focused purely on charge scheduling such as \citet{abdelwahed_evaluating_2020} and \citet{liu_robust_2022} are not applicable under all conditions, since they assume buses run on schedule and charging can take place within scheduled layover time. There is a need for models that can provide useful output even under challenging conditions when avoiding departure delays is impossible, so that bus operators are provided with actionable instructions rather than being told a problem instance is infeasible. Such models should be accompanied with efficient algorithms so they can be run repeatedly as conditions evolve or forecasts of ridership, traffic, and weather conditions improve.

Based on this research gap, we propose one such model for charging scheduling. Our approach emphasizes quality of service by seeking to minimize departure delays across all trips. We precisely quantify the impact of queuing at chargers and ensure delays are propagated across trips. The resulting model and methods for its solution are presented next in Sections \ref{sec:methodology}--\ref{sec:heuristic}.

\section{Mathematical Programming Formulation}
\label{sec:methodology}
\subsection{Problem Setting, Assumptions, and Modeling Approach}
\begin{figure}[ht]
    \centering
    \begin{tikzpicture}
        \Vertex[label=$l_1$, x=0,y=4, size=0.5, color=white]{l1}
        \Vertex[label=$l_2$, x=6,y=4, size=0.5, color=white]{l2}
        
        \Edge[Direct, loopposition=180, loopsize=3cm, color=orange](l1)(l1);
        \node[orange] at (-1, 3) {\large Route B};
        \node[rotate=90] at (-2.1, 4) {\includegraphics[width=1cm]{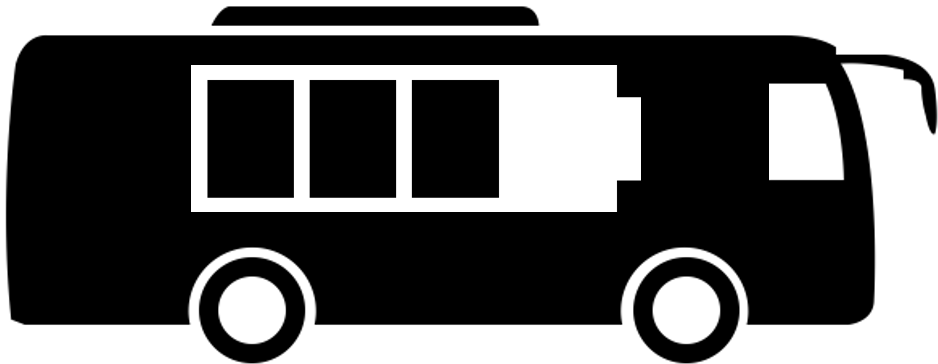}};
        \Edge[Direct, bend=15, color=tlblue](l1)(l2)
        \Edge[Direct, bend=15, color=tlblue](l2)(l1)
        \node[tlblue] at (3, 4) {\large Route A};
        \node at (3, 4.75) {\includegraphics[width=1cm]{images/beb_icon_right.png}};
        \node[rotate=180] at (3, 3.25) {\includegraphics[width=1cm]{images/beb_icon_right.png}};
        
        \node[tlgreen] at (0.1, 4.7) {\Huge \Lightning};
        \node[tlgreen] at (6.1, 4.7) {\Huge \Lightning};
        
    \end{tikzpicture}
    \caption{Simple example of the type of BEB system considered in this study.}
    \label{fig:system_sketch}
\end{figure}

We consider a general setting in which fast chargers with predetermined power outputs have already been installed at some terminals of a BEB system. Figure \ref{fig:system_sketch} shows a simple illustration of the type of transit network considered in this study. In this basic example, two bus routes (A and B) operate across two terminals and each terminal has a single charger installed. A bus may recharge at either of the two charging sites when out of service in between trips, as long as the bus is at the corresponding terminal and the charger is not already occupied. Note that in our approach, a single bus may serve more than one route (referred to as \textit{interlining}) and/or use multiple different chargers during a day. A single charger also may serve buses on any number of different routes, so that infrastructure is shared as efficiently as possible. Bus trips may take place on a one-way route between distinct terminals (i.e., Route A) or a loop route that starts and ends at the same terminal (i.e., Route B).

Our approach relies on the following assumptions:

\begin{enumerate}
    \item Buses rely on slow charging overnight at the depot in addition to fast charging during the day. Each bus has a known state of charge when entering service at the start of the day and a minimum state of charge it must have when returning to the depot at the end of service. \label{assume:depot}
    \item Trip schedules and vehicle blocks have been determined in advance and cannot be altered. The agency's primary goal is to adhere to the posted schedule as much as possible, i.e., to minimize delays. \label{assume:fixedsched}
    \item Buses may only charge at terminals when passengers are not onboard, never at intermediate stops. \label{assume:nointermed}
    \item Deadheading to chargers is not allowed -- either a bus is close enough to a charger to use it at the end of a trip, or it cannot charge at all. \label{assume:dh}
    \item Charger locations are established in advanced. At most one charger is located at each terminal. \label{assume:singlecharger}
    \item All buses in the network have sufficient battery capacity to complete all trips in between successive charging opportunities. \label{assume:batteries}
    \item Charging behavior is linear, proportional to the maximum power output of each charger. \label{assume:linearchg}
\end{enumerate}

Assumption \ref{assume:depot} acknowledges the different technologies, time scales, and constraints faced by transit agencies during depot versus opportunity charging. Agencies can be expected to have different priorities for charging scheduling in these different environments. When buses use opportunity chargers during the service day, maintaining schedule adherence is the greatest priority, whereas overnight charging should allow more time to optimize for other objectives such as energy costs and battery health. Bus bases are also less likely to be subject to queuing for chargers due to the much lower cost of lower-power chargers.

Assumption \ref{assume:fixedsched} reflects our intended problem setting as described in Section \ref{sec:introduction}. Assumption \ref{assume:nointermed} ensures that the current passenger level of service is maintained. Assumption \ref{assume:dh} is common in the recharging literature and simplifies the model by making cumulative energy consumption independent of any charging decisions and their corresponding binary decision variables. Assumption \ref{assume:singlecharger} is likely accurate for many agencies and helps simplify the formulation of our mathematical program as well as the heuristic method developed to solve it. In Section \ref{sec:conclusion}, we discuss how this assumption can be relaxed without much computational penalty to apply our methods to networks where buses must choose between multiple available chargers at a given time.

The last two assumptions concern battery range and charging behavior. Assumption \ref{assume:batteries} is necessary to ensure operations are always feasible with fast chargers; any such infeasible blocks can be discarded from the analysis. Assumption \ref{assume:linearchg} is another common assumption made to limit computational complexity. It can provide a good enough approximation even if the real charging behavior is nonlinear, since charging is typically linear until the state of charge exceeds about 80\% when power must be reduced to protect the battery \citep{Montoya2017}.

Following these assumptions, suppose we are given a set of BEBs $\setV$, each of which is scheduled to complete a specific set of trips; that is, every vehicle has been assigned to a predetermined block. We represent each trip as a tuple $i = (b_i, n_i)$ where $b_i$ is the bus ID and $n_i$ is the trip number; i.e., bus $A$ completes trips $(A, 1), (A, 2), (A, 3)$, etc. Let $\setT$ be the set of all trips completed by any bus. We will use the shorthand indices $i$ and $j \in \setT$ throughout this formulation to limit indices, but keep in mind that each of these refer to a specific trip completed by a specific bus as part of a predetermined block. For each trip $i \in \setT$, the timetable provides a scheduled departure time $\sigma_i$ and scheduled duration $\tau_i$ (so that the scheduled end time is $\sigma_i + \tau_i$). If a trip $i$ has an immediate predecessor in its block, that predecessor is denoted $i^-$. Let $\theta(i)$ be the set of trips completed by the same bus prior to $i$ (i.e., all of its predecessor trips).

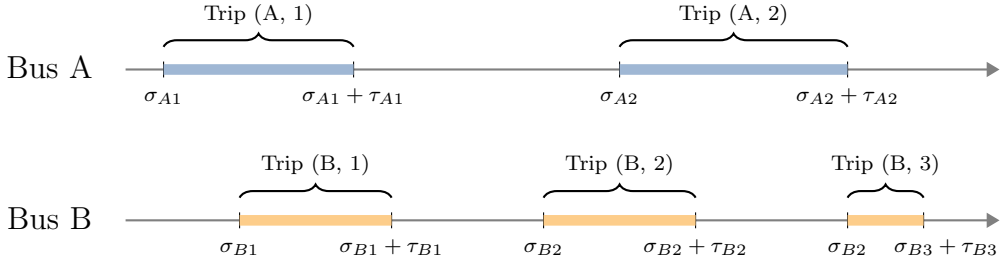
\begin{figure}[ht]
    \centering
    \begin{tikzpicture}
    \draw[gray, thick, -Triangle] (-0.5,0) -- (\ImageWidth,0) node[font=\scriptsize,below left=3pt and -8pt]{};
    \node[font=\normalsize, text height=1.75ex, text depth=.5ex] at (-1.5,0) {Bus B};
    
    \foreach \x/\descr in {1/\sigma_{B1},3/\sigma_{B1}+\tau_{B1},5/\sigma_{B2},7/\sigma_{B2}+\tau_{B2},9/\sigma_{B2},10/\qquad \sigma_{B3}+\tau_{B3}}
    {\draw (\x cm,3pt) -- (\x cm,-3pt);
    \node[font=\scriptsize, text height=1.75ex,
    text depth=.5ex] at (\x,-.3) {$\descr$};}

    \draw[busbfaded, line width=4pt] (1, 0) -- (3, 0);
    \draw [thick,decorate,decoration={brace,amplitude=5pt}] (1,0.3)  -- +(2,0) 
           node [black,midway,above=4pt, font=\scriptsize] {Trip (B, 1)};

    \draw[busbfaded, line width=4pt] (5, 0) -- +(2, 0);
    \draw[thick,decorate,decoration={brace,amplitude=5pt}] (5,0.3)  -- +(2,0) 
           node [black,midway,above=4pt, font=\scriptsize] {Trip (B, 2)};

    \draw[busbfaded, line width=4pt] (9, 0) -- +(1, 0);
    \draw[thick,decorate,decoration={brace,amplitude=5pt}] (9,0.3)  -- +(1,0) 
           node [black,midway,above=4pt, font=\scriptsize] {Trip (B, 3)};
    
    \draw[gray, thick, -Triangle] (-0.5,2) -- (\ImageWidth,2) node[font=\scriptsize,below left=3pt and -8pt]{};
    \node[font=\normalsize, text height=1.75ex, text depth=.5ex] at (-1.5,2) {Bus A};

    \foreach \x/\descr in {0/\sigma_{A1},2.5/\sigma_{A1}+\tau_{A1},6/\sigma_{A2},9/\sigma_{A2}+\tau_{A2}}
    {\draw ([yshift=2cm]\x cm, 3pt) -- ([yshift=2cm]\x cm,-3pt);
    \node[font=\scriptsize, text height=1.75ex,
    text depth=.5ex] at ([yshift=2cm] \x,-.3) {$\descr$};};

    \draw[busafaded, line width=4pt] ([yshift=2cm] 0, 0) -- +(2.5, 0);
    \draw [thick,decorate,decoration={brace,amplitude=5pt}] ([yshift=2cm] 0,0.3)  -- +(2.5,0) 
           node [black,midway,above=4pt, font=\scriptsize] {Trip (A, 1)};

    \draw[busafaded, line width=4pt] ([yshift=2cm] 6, 0) -- +(3, 0);
    \draw[thick,decorate,decoration={brace,amplitude=5pt}] ([yshift=2cm] 6,0.3)  -- +(3,0) 
           node [black,midway,above=4pt, font=\scriptsize] {Trip (A, 2)};
    
    \end{tikzpicture}
    \caption{Timeline of operations for two example buses.}
    \label{fig:param_timeline}
\end{figure}

Figure \ref{fig:param_timeline} illustrates a simple example of how to interpret these parameters on a timeline for two buses labeled A and B. The colored segments of the timeline represent times when the buses are scheduled to be in service and unavailable to charge. Considering the trips served by bus B in this example, we have $\theta((B, 3)) = \{(B, 1), (B, 2)\}$, $\theta((B, 2)) = \{(B, 1)\}$, and $\theta((B,1)) = \emptyset$. Likewise, $(B, 3)^- = (B, 2)$ and $(B, 2)^- = (B, 1)$.

Our approach to charging scheduling is aimed at accurately tracking trip departure delays. Tracking and propagating delay within a mixed-integer linear programming paradigm is challenging because of the underlying queuing behavior. To address this challenge, we define multiple continuous variables for each trip that track the time of certain events. This indexing strategy can be seen as similar to the ``discrete event optimization'' approach of \citet{abdelwahed_evaluating_2020}. For each trip $i$, let $d_i$ be a decision variable giving its departure delay and $p_i$ be a decision variable giving the plug-in time after completion of trip $i$, i.e., the time the bus that completes $i$ connects to a charger. To ensure that delay is tracked properly, we also require that if there is no charging after $i$, the value of $p_i$ equals the completion time of trip $i$, i.e., $ p_i = d_i + \sigma_i + \tau_i$. Let $t^l_i$ be the amount of time spent charging after trip $i$ at charger $l$.

\begin{figure}[ht]
    \centering
    \begin{tikzpicture}
    \draw[gray, thick, -Triangle] (-0.5,0) -- (\ImageWidth,0) node[font=\scriptsize,below left=3pt and -8pt]{};
    \node[font=\normalsize, text height=1.75ex, text depth=.5ex] at (-1.5,0) {Bus B};

    \foreach \x/\descr in {3/p_{B1},4.5/p_{B1}+t^{l_1}_{B1}}
    {\draw (\x cm, 3pt) -- (\x cm,-3pt);
    \node[font=\scriptsize, text height=1.75ex,
    text depth=.5ex] at (\x,.4) {$\descr$};};

    \draw[busbfaded, line width=4pt] (1, 0) -- (3, 0);
    \draw[pastelgreen, line width=4pt] (3, 0) -- +(1.5, 0);

    \draw[busbfaded, line width=4pt] (5, 0) -- +(2, 0);

    \draw[busbfaded, line width=4pt] (9, 0) -- +(1, 0);

    \foreach \x/\descr in {1/\sigma_{B1},3/\sigma_{B1}+\tau_{B1},5/\sigma_{B2},7/\sigma_{B2}+\tau_{B2},9/\sigma_{B2},10/\qquad \sigma_{B3}+\tau_{B3}}
    {\draw (\x cm,3pt) -- (\x cm,-3pt);
    \node[font=\scriptsize, text height=1.75ex,
    text depth=.5ex] at (\x,-.3) {$\descr$};}
    
    \draw[gray, thick, -Triangle] (-0.5,2) -- (\ImageWidth,2) node[font=\scriptsize,below left=3pt and -8pt]{};
    \node[font=\normalsize, text height=1.75ex, text depth=.5ex] at (-1.5,2) {Bus A};

    \foreach \x/\descr in {4.5/p_{A1},7/p_{A1}+t^{l_1}_{A1}}
    {\draw ([yshift=2cm]\x cm, 3pt) -- ([yshift=2cm]\x cm,-3pt);
    \node[font=\scriptsize, text height=1.75ex,
    text depth=.5ex] at ([yshift=2cm]\x,.3) {$\descr$};};

    \draw[busafaded, line width=4pt] ([yshift=2cm] 0, 0) -- +(2.5, 0);
    \draw[pastelgreen, line width=4pt] (4.5, 2) -- +(2.5, 0);
    \draw [black,thick,decorate,decoration={brace,amplitude=5pt}] ([yshift=2cm] 6,0.5)  -- +(1,0) 
           node [black,midway,above=4pt, font=\scriptsize] {$d_{A2}$};

    \draw[busafaded, line width=4pt] ([yshift=2cm] 7, 0) -- +(3, 0);

    \foreach \x/\descr in {0/\sigma_{A1},2.5/\sigma_{A1}+\tau_{A1},6/\sigma_{A2},9/\sigma_{A2}+\tau_{A2}}
    {\draw ([yshift=2cm]\x cm, 3pt) -- ([yshift=2cm]\x cm,-3pt);
    \node[font=\scriptsize, text height=1.75ex,
    text depth=.5ex] at ([yshift=2cm] \x,-.3) {$\descr$};};

    \end{tikzpicture}
    \caption{Illustration of variable meanings for $d_i$, $t^l_i$, and $p_i$ with their relation to schedule parameters $\sigma_i$ and $\tau_i$.}
    \label{fig:var_timeline}
\end{figure}
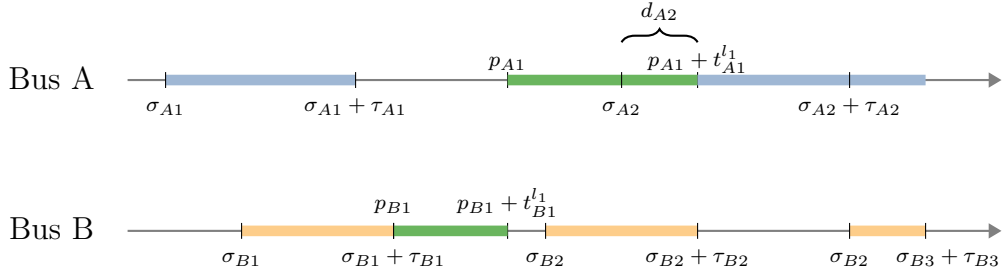

To illustrate exactly how these variables are related to one another and the problem parameters, Figure \ref{fig:var_timeline} shows a simple example based on the schedule of Figure \ref{fig:param_timeline}. Suppose buses A and B both use the same charger $l_1$ after completing their first trip of the day. Bus B arrives first, immediately plugs in at time $p_{B1}=\sigma_{B1}+\tau_{B1}$ and charges for time $t^{l_1}_{B1}$. It then starts trip $(B, 2)$ on time, since $\sigma_{B1}+\tau_{B1}+t^{l_1}_{B1}$ is earlier than the scheduled start time $\sigma_{B2}$ of its next trip. When Bus B finishes charging, the charger becomes available and Bus A plugs in at time $p_{A1} = p_{B1}+t^{l_1}_{B1}$. Bus A charges for length $t^{l_1}_{A1}$. Figure \ref{fig:var_timeline} shows that Bus A departs its second trip late because of the time it spent waiting at the charger; specifically, its departure delay is $d_{A2} = p_{A1} + t^{l_1}_{A1} - \sigma_{A2}$.

Table \ref{tab:notation} compiles all the set, parameter, and decision variable definitions used in this work. Sections \ref{sec:obj}--\ref{sec:formulation} next describe the formulation of the mixed-integer linear program we developed for recharging scheduling, including the objective function and the constraints that track queuing, delays, and battery levels.

\begin{table}[!ht]
	\begin{center}
		\caption{Notation definitions.}
		\begin{tabular}{|c|l|}
			\hline 
			\multicolumn{2}{|l|}{\textbf{\textit{Sets}}} \\ \hline
			$\setV$ & Buses \\
            $\setT$ & Trips \\
            $\setC$ & Chargers \\
            $\setA$ & Arcs \\
            \hline
			\multicolumn{2}{|l|}{\textit{\textbf{Decision Variables}}} \\ \hline
			$x^l_i$ & Binary variable indicating whether the bus on trip $i$ uses charger $l$ afterwards \\
			$y^l_{ij}$ & Binary variable indicating whether charging arc $(l, i, j)$ is used \\
            $d_i$ & Delay at start of trip $t$ \\
            $p_i$ & Plugin time of charging after trip $i$ \\
			$t^l_i$ & Charging duration at $l$ after completion of trip $i$ \\ \hline
			\multicolumn{2}{|l|}{\textit{\textbf{Parameters}}} \\ \hline
            $b_i$ & ID of bus that completes trip $i$ \\
            $n_i$ & Trip number of trip $i$ \\
			$\alpha_i$ & Minimum total charge gain by end of trip $i$ \\
            $\beta_i$ & Maximum total charge gain by end of trip $i$ \\
			$\rho^l$ & Power output of charger $l$ \\
			$\overline{t^l_i}$ & Maximum charging duration at site $l$ after trip $i$ \\
			$\sigma_{i}$ & Scheduled start time of trip $i$ \\
			$\tau_{i}$ & Scheduled duration of trip $i$ \\
            $M$ & ``Big M'', arbitrary large number \\
			\hline
		\end{tabular}
		\label{tab:notation}
	\end{center}
\end{table}

\subsection{Objective Function}
\label{sec:obj}
\begin{equation}
    \min_{\substack{d_i, p_i, t^l_i, \\x^l_i, y^l_{ij}}} \quad \sum_{i \in \setT} d_i \label{eq:mipobj}
\end{equation}

The objective \eqref{eq:mipobj} is to minimize the total amount of departure delay across all trips $i \in \setT$ served by BEBs. While the objective is straightforward, setting the appropriate value of $d_i$ for each trip $i$ requires carefully formulated constraints, as described in the following subsections.

\subsection{Queue Tracking and Delay Propagation}

Our queuing and delay tracking constraints are based on the relationships between plugin time, charging time, and delay depicted in Figure \ref{fig:var_timeline}. We first describe how to correctly set the plugin time. Recall that if bus $b_i$ does not charge after trip $i$, we define the plugin time $p_i$ to be the end time of trip $i$. Since charging can only increase the value of $p_i$, the trip completion time provides a universal lower bound on $p_i$, encoded by Constraints \eqref{eq:pluginlb}.

\begin{align}
    p_i &\geq d_i + \sigma_i + \tau_i \qquad & \forall \enskip i \in \setT \label{eq:pluginlb}
\end{align}


If $b_i$ does use some charger $l$ after trip $i$, then setting its plugin time correctly is more complex because it may have to queue. We therefore introduce some auxiliary binary variables and related constraints. Let $\setC$ be the set of all chargers that are installed. We say that a charger $l \in \setC$ \textit{serves} a trip $i \in \setT$ if $b_i$ uses $l$ after finishing trip $i$. Let $y^l_{ij}$ be a binary decision variable that equals 1 if charger $l$ serves trip $j$ immediately after trip $i$ and 0 otherwise. By ``immediately,'' we mean only that no other trips are served in between these two; there may be a large gap between the charging completion time for trip $i$ and the plugin time of trip $j$. Likewise, let $x^l_i = 1$ if charger $l$ serves trip $i$ and 0 otherwise. Then we account for queuing at the charger with Constraints \eqref{eq:pluginbigm}:

\begin{align}
    p_j &\geq p_i + t^l_i - M(1 - y^l_{ij}) \quad & \forall \enskip (l, i, j) \in \setA \label{eq:pluginbigm}
\end{align}

Constraints \eqref{eq:pluginbigm} impose additional lower bounds on plugin time to account for queuing. These ``big M'' constraints encode the logic that if trip $j$ is served by charger $l$ immediately after trip $i$, then $b_j$ cannot plug in until $b_i$ has unplugged. Written as an equivalent logical constraint, we would have

\begin{align}
    & y^l_{ij} = 1 \implies p_j \geq p_i + t_i \quad & \forall \enskip (l, i, j) \in \setA \label{eq:pluginlogic}
\end{align}

In Constraints \eqref{eq:pluginbigm} and \eqref{eq:pluginlogic}, $\setA$ is the set of ``arcs'' connecting trips $i$ and $j$ that may be served sequentially by charger $l$. These arcs are discussed in more detail in Section \ref{sec:chgseq}.

With the plugin time values constrained properly, it is straightforward to set the delay for each trip. Note that whether or not bus $b_i$ charges after $i$, it is ready to re-enter service at time $p_i + \max_{l \in \setC} \left\{t^l_i\right\}$. The departure delay of trip $i$ is then $d_i = \max \left(0, p_{i^-} + \max_{l \in \setC} \left\{t^l_{i^-}\right\} - \sigma_i \right)$. This maximum operator can be represented exactly with Constraints \eqref{eq:delay} and \eqref{eq:delaynonneg}.

\begin{align}
    d_{i} &\geq p_{i^-} + t^l_{i^-} - \sigma_{i} & \forall \enskip l \in \setC, \enskip i \in \setT: n_i \geq 2 \label{eq:delay} \\
    d_i &\geq 0 & \forall \enskip i \in \setT \label{eq:delaynonneg}
\end{align}

\subsection{Charger Sequencing}
\label{sec:chgseq}
The binary variables $y^l_{ij}$ explicitly track the sequence of trips that are served by each charger in the network. Our model therefore includes constraints to ensure that the optimal values of $y^l_{ij}$ encode a valid sequence connecting all trips served by each charger. By a valid sequence, we mean that (1) any time a charger is used (i.e., $x^l_i=1$), that trip appears somewhere in the sequence, and (2) every trip $i$ served by charger $l$ has exactly one trip before and one trip after it in the sequence.

To formulate these constraints, we model the charging sequence for each charger $l$ as a path through a network in which each node corresponds to a trip (a charging opportunity) and arcs represent feasible charging connections. That is, if trip $j$ can be served by the charger immediately after trip $i$, then arc $(l, i, j)$ is included in the network. In general, a directed arc joins every pair of edges for the same bus (since we can't travel back in time to charge at, say, trip 1 after trip 3). On the other hand, an undirected edge joins any pair of trips completed by different buses, because we don't know \textit{a priori} which order these trips will be optimally served in. For instance, it could be optimal to use arc $(l, i, j)$ in the solution even if trip $j$ is scheduled to end before trip $i$, because trip $j$ is delayed in the optimal solution.

We also introduce two dummy nodes to model the initial (node $s$) and final (node $t$) idle state of each charger. These nodes are necessary because the first and last trips to be served by each charger have no predecessor and successor trip, respectively, so they require slightly different constraints. Moreover, there is no way to know \textit{a priori} which trips will be first and last, so we instead construct dummy nodes to handle these special cases. Now, a feasible sequence for charger $l$ corresponds to a path through this virtual network from $s$ to $t$. Note that this network model is similar to and inspired by those used in vehicle scheduling approaches (e.g., the maximum flow formulation of the VSP \citep{Ceder2007}).

\begin{figure}[ht]
    \centering
    \begin{subfigure}{.45\textwidth}
      \centering
      \begin{tikzpicture}[scale=.8]
        \Vertex[color=gray!40,label=s, x=-0.25, y=1.5]{s}
        \Vertex[color=gray!40,label=t, x=6.25, y=1.5]{t}
    
        \Vertex[color=busafaded,label=\textcolor{black}{A1},x=2,y=3]{A1}
        \Vertex[color=busafaded,label=\textcolor{black}{A2},x=4,y=3]{A2}
        \Vertex[color=busbfaded,label=\textcolor{black}{B1},x=1,y=0]{B1}
        \Vertex[color=busbfaded,label=\textcolor{black}{B2},x=3,y=0]{B2}
        \Vertex[color=busbfaded,label=\textcolor{black}{B3},x=5,y=0]{B3}

        \Edge[Direct, style=dashed, color=gray](s)(A1)
        \Edge[Direct, style=dashed, color=gray](s)(A2)
        \Edge[Direct, style=dashed, color=gray](s)(B1)
        \Edge[Direct, style=dashed, color=gray](s)(B2)
        \Edge[Direct, style=dashed, color=gray](s)(B3)
        \Edge[Direct, style=dashed, color=gray](A1)(t)
        \Edge[Direct, style=dashed, color=gray](A2)(t)
        \Edge[Direct, style=dashed, color=gray](B1)(t)
        \Edge[Direct, style=dashed, color=gray](B2)(t)
        \Edge[Direct, style=dashed, color=gray](B3)(t)

        \Edge[Direct](A1)(A2)
        \Edge[Direct](B1)(B2)
        \Edge[Direct](B2)(B3)
        \Edge[Direct, bend=-30](B1)(B3)

        \Edge[](A1)(B1)
        \Edge[](A1)(B2)
        \Edge[](A1)(B3)
        \Edge[](A2)(B1)
        \Edge[](A2)(B2)
        \Edge[](A2)(B3)
        
    \end{tikzpicture}
    \caption{Network structure.}
    \label{fig:spnetwork}
    \end{subfigure}
    \begin{subfigure}{.45\textwidth}
    \centering
        \begin{tikzpicture}[scale=.8]
        \Vertex[color=gray!40,label=s, x=-0.25, y=1.5]{s}
        \Vertex[color=gray!40,label=t, x=6.25, y=1.5]{t}
    
        \Vertex[color=busafaded,label=\textcolor{black}{A1},x=2,y=3]{A1}
        \Vertex[color=busafaded,label=\textcolor{black}{A2},x=4,y=3]{A2}
        \Vertex[color=busbfaded,label=\textcolor{black}{B1},x=1,y=0]{B1}
        \Vertex[color=busbfaded,label=\textcolor{black}{B2},x=3,y=0]{B2}
        \Vertex[color=busbfaded,label=\textcolor{black}{B3},x=5,y=0]{B3}

        \Edge[Direct, style=dashed, color=gray](s)(A1)
        \Edge[Direct, style=dashed, color=gray](s)(A2)
        \Edge[Direct, color=red](s)(B1)
        \Edge[Direct, style=dashed, color=gray](s)(B2)
        \Edge[Direct, style=dashed, color=gray](s)(B3)
        \Edge[Direct, color=red](A1)(t)
        \Edge[Direct, style=dashed, color=gray](A2)(t)
        \Edge[Direct, style=dashed, color=gray](B1)(t)
        \Edge[Direct, style=dashed, color=gray](B2)(t)
        \Edge[Direct, style=dashed, color=gray](B3)(t)

        \Edge[Direct](A1)(A2)
        \Edge[Direct](B1)(B2)
        \Edge[Direct](B2)(B3)
        \Edge[Direct, bend=-30](B1)(B3)

        \Edge[Direct, color=red](B1)(A1)
        \Edge[](A1)(B2)
        \Edge[](A1)(B3)
        \Edge[](A2)(B1)
        \Edge[](A2)(B2)
        \Edge[](A2)(B3)
        
    \end{tikzpicture}
    \caption{Example feasible solution from Figure \ref{fig:var_timeline}}
    \label{fig:spexample}
    \end{subfigure}
    
    \caption{Example network for charger $l_1$ serving buses A and B.}
    \label{fig:fullnetworkfig}
    
\end{figure}
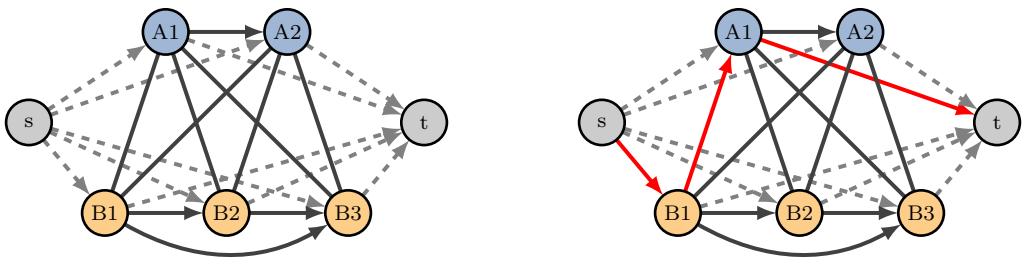

Figure \ref{fig:fullnetworkfig} illustrates our network model for the simple example of two buses and one charger previously depicted in Figures \ref{fig:param_timeline} and \ref{fig:var_timeline}. Figure \ref{fig:spnetwork} displays the network structure including the trip and dummy nodes as well as the arc set $\setA$ (where arrows indicate directed arcs and lines indicate undirected edges). Figure \ref{fig:spexample} shows in red the path corresponding to the charging decisions from Figure \ref{fig:var_timeline}, where buses charge at $l$ after trips $(B, 1)$ and $(A, 1)$. The corresponding arc variable values are $y^l_{sB1} = y^l_{B1A1} = y^l_{A1t} = 1$, and $y^l_{ij} = 0$ for all other arcs.

Following this network representation, a feasible sequence for each charger can be enforced using the familiar network flow constraints \eqref{eq:leavedummy}--\eqref{eq:flowbalance}.

\begin{align}
    & \sum_{i \in \setT} y^l_{si} = 1 & \forall \enskip l \in \setC \label{eq:leavedummy} \\
    & \sum_{i \in \setT} y^l_{it} = 1 & \forall \enskip l \in \setC \label{eq:returndummy} \\
    & \sum_{j: (l, j, i) \in \setA} y^l_{ji} - \sum_{j: (l, i, j) \in \setA} y^l_{i j} = 0 & \forall \enskip l \in \setC, \forall \enskip i \in \setT \label{eq:flowbalance}
\end{align}

Constraints \eqref{eq:leavedummy} and \eqref{eq:returndummy} ensure that exactly one arc leaves the source dummy node and exactly one arc arrives at the sink dummy node, respectively, for each charger. Constraints \eqref{eq:flowbalance} ensure connectivity of the charging sequence; the number of arcs entering and leaving each trip node must be equal. Notably, the network flow constraints \eqref{eq:leavedummy}--\eqref{eq:flowbalance} do not preclude the existence of subtours in paths through the network. It is not necessary to include subtour elimination constraints because any solution containing a subtour cannot possibly be optimal for our problem; in fact, it would result in unbounded delay due to Constraints \eqref{eq:pluginbigm} and \eqref{eq:delay}.

Finally, we need to establish the relationship between $x^l_i$ and $y^l_{ij}$ variables. Charging after trip $i$ corresponds to visiting its trip node in Figure \ref{fig:fullnetworkfig}, meaning one arc leaves (or, equivalently, enters) that node, so we have:

\begin{align}
    & x^l_i = \sum_{j: (l, i, j) \in \setA} y^l_{ij} & \forall \enskip l \in \setC, \enskip i \in \setT
\end{align}

\subsection{State of Charge Management}
\label{sec:chgconstr}
The final subsystem of constraints tracks charging throughout the service day to make sure that each bus's state of charge stays within an acceptable range (e.g.,  between 20\% and 95\%). First, we need to relate the continuous charging duration $t^l_i$ to the binary $x^l_i$. $t^l_i$ must be forced to zero if $x^l_i = 0$; otherwise, it can take any feasible value. Let $\overline{t^l_i}$ be the maximum possible length of charging at charger $l$ after trip $i$. Let $u^{\text{max}}_i$ be the usable battery capacity of the bus serving trip $i$ and $\rho^l$ be the power output of charger $l$. We set the value of $\overline{t^l_i}$ with Equation \eqref{eq:tmax}:

\begin{align}
    & \overline{t^l_i} = \begin{cases} u^{\text{max}}_i / \rho^l & \text{if trip $i$ ends at $l$}  \\ 0 & \text{otherwise}  \end{cases} & \forall \enskip l \in \setC, \enskip i \in \setT \label{eq:tmax}
\end{align}

Equation \eqref{eq:tmax} reflects that bus $b_i$ can only use charger $l$ if trip $i$ ends at $l$. Additionally, the maximum possible charging duration is the battery capacity of $b_i$ divided by the charging rate of $l$. With $\overline{t^l_i}$ defined as such, $t^l_i$ is restricted by Constraints \eqref{eq:chgtime}:

\begin{align}
    & 0 \leq t^l_{i} \leq \overline{t^l_{i}} x^l_i & \forall \enskip l \in \setC, \enskip i \in \setT \label{eq:chgtime}
\end{align}

We also need constraints to ensure that each bus charges enough during the day in order to complete its scheduled service. Following Assumption \ref{assume:linearchg}, with linear charging it is straightforward to determine lower and upper bounds $(\alpha_i, \beta_i)$ on the cumulative charge gain through the end of each trip $i$. The complete state-of-charge constraints are given by Constraints \eqref{eq:soc}.

\begin{align}
    & \alpha_i \leq \sum_{j \in \theta(i)} \sum_{l \in \setC} \rho^l t^l_j \leq \beta_i & \forall \enskip i \in \setT \label{eq:soc}
\end{align}

To understand how $\alpha_i$ and $\beta_i$ can be calculated for all trips, consider a simple example. Suppose bus A must keep its battery level between 50 and 200 kWh all day and enters service with 125 kWh of charge. Suppose its first two trips each require 50 kWh of energy. Then, its charge level when starting trip 2 must be no more than its maximum (200 kWh) and no less than enough to complete trip 2 without falling below its minimum ($50 + 50 = 100$ kWh). Since its battery level is $125 - 50 = 75$ kWh after completing trip 1, it must gain between 25 and 125 kWh of charge after that trip. Accordingly, $l_{A1} = 25$ and $u_{A1} = 125$. This simple example for trip (A, 1) easily generalizes across all trips.

\subsection{Complete MILP Formulation}
\label{sec:formulation}
The complete mixed-integer linear programming formulation for recharging scheduling as developed in Sections \ref{sec:obj}--\ref{sec:chgconstr} is given by Equations \eqref{eq:mipobj2}--\eqref{eq:mipend}:

\begin{align}
    \min_{\substack{d_i, p_i, t^l_i, \\x^l_i, y^l_{ij}}} \quad & \sum_{i \in \setT} d_i & \label{eq:mipobj2} \\
    \text{s.t.} \quad & d_{i} \geq p_{i^-} + t^l_{i^-} - \sigma_{i} & \forall \enskip l \in \setC, \enskip i \in \setT: n_i \geq 2 \\
    & p_{i} \geq \sigma_{i} + d_{i} + \tau_{i} & \forall \enskip i \in \setT \label{eq:pluginlbfull} \\
    & p_{j} \geq p_{i} + t^l_{i} - M \left(1 - y^l_{i j} \right) & \forall \enskip (l, i, j) \in \setA \label{eq:bigmfull} \\
    & \sum_{i \in \setT} y^l_{si} = 1 & \forall \enskip l \in \setC \\
    & \sum_{i \in \setT} y^l_{it} = 1 & \forall \enskip l \in \setC \\
    & \sum_{j: (l, j, i) \in \setA} y^l_{ji} - \sum_{j: (l, i, j) \in \setA} y^l_{i j} = 0 & \forall \enskip l \in \setC, \enskip i \in \setT \\
    & x^l_i = \sum_{j: (l, i, j) \in \setA} y^l_{ij} & \forall \enskip l \in \setC, \enskip i \in \setT \\
    & 0 \leq t^l_{i} \leq \overline{t^l_{i}} x^l_i & \forall \enskip l \in \setC, \enskip i \in \setT \label{eq:chgtimefull} \\
    & \alpha_i \leq \sum_{j \in \theta(i)} \sum_{l \in \setC} \rho^l t^l_j \leq \beta_i & \forall \enskip i \in \setT \label{eq:socfull} \\
    & d_i \geq 0 & \forall \enskip i \in \setT \\
    & x^l_i \in \{0, 1\} & \forall \enskip l \in \setC, \enskip i \in \setT \\
    & y^l_{ij} \in \{0, 1\} & \forall \enskip (l, i, j) \in \mathcal{A} \label{eq:mipend}
\end{align}

\section{Exact Solution Method: Combinatorial Benders Decomposition}
\label{sec:benders}
The formulation \eqref{eq:mipobj2}--\eqref{eq:mipend} is difficult to solve for large instances, largely due to the binary variables $y^l_{ij}$ and the corresponding big-M constraints \eqref{eq:bigmfull}. Because the number of $y$ variables and big-M constraints scales with the square of the number of trips, solving the problem with an off-the-shelf MIP solver is not possible for many large instances. We develop two strategies to deal with this computational challenge. In this section, we describe an exact solution approach based on Combinatorial Benders (CB) decomposition, which circumvents the typical issues of weak linear programming relaxations caused by big-M constraints. Section \ref{sec:heuristic} later describes a polynomial-time heuristic with randomization that can generate a large number of good solutions quickly. The heuristic helps accelerate the convergence of the CB algorithm on smaller problems and also shows good performance as a standalone method for difficult problems on real networks where exact algorithms are unacceptably slow.

\subsection{Overview of Combinatorial Benders Decomposition}
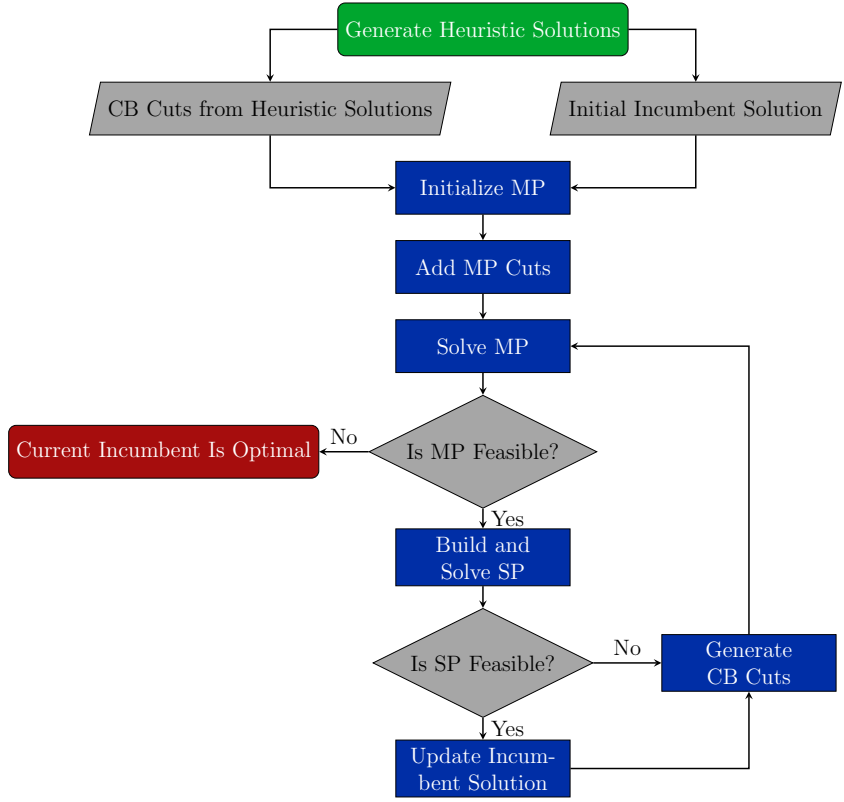
\begin{figure}[ht]
    \begin{center}
    \scalebox{0.7}{
        \begin{tikzpicture}[node distance=2cm]
        \node (start) [start] {Generate Heuristic Solutions};
        \node (cbheur) [io, below of=start, yshift=0.5cm, xshift=-4cm] {CB Cuts from Heuristic Solutions};
        \node (inc0) [io, below of=start, yshift=0.5cm, xshift=4cm] {Initial Incumbent Solution};
        \node (mp0) [process, below of=start, yshift=-1cm] {Initialize MP};
        \node (mpcuts) [process, below of=mp0, yshift=0.5cm] {Add MP Cuts};
        \node (mpsol) [process, below of=mpcuts, yshift=0.5cm] {Solve MP};
        \node (mpfeas) [decision, below of=mpsol] {Is MP Feasible?};
        \node (stop) [stop, left of=mpfeas, xshift=-4cm] {Current Incumbent Is Optimal};
        \node (spsol) [process, below of=mpfeas] {Build and Solve SP};
        \node (spfeas) [decision, below of=spsol] {Is SP Feasible?};
        \node (incupd) [process, below of=spfeas] {Update Incumbent Solution};
        \node (cbgen) [process, right of=spfeas, xshift=3cm] {Generate CB Cuts};
        
        \draw [arrow] (start) -| (cbheur);
        \draw [arrow] (start) -| (inc0);
        \draw [arrow] (cbheur) |- (mp0);
        \draw [arrow] (inc0) |- (mp0);
        \draw [arrow] (mp0) -- (mpcuts);
        \draw [arrow] (mpcuts) -- (mpsol);
        \draw [arrow] (mpsol) -- (mpfeas);
        \draw [arrow] (mpfeas) -- node[anchor=south] {No} (stop);
        \draw [arrow] (mpfeas) -- node[anchor=west] {Yes} (spsol);
        \draw [arrow] (spsol) -- (spfeas);
        \draw [arrow] (spfeas) -- node[anchor=west] {Yes} (incupd);
        \draw [arrow] (spfeas) -- node[anchor=south] {No} (cbgen);
        \draw [arrow] (incupd) -| (cbgen);
        \draw [arrow] (cbgen) |- (mpsol);
        \end{tikzpicture}
    }
    \end{center}
    \caption{Overview of Combinatorial Benders solution process.}
    \label{fig:cbflowchart}
\end{figure}

Figure \ref{fig:cbflowchart} shows a high-level overview of the CB algorithm we implemented. We begin by using a randomized heuristic (described in detail in Section \ref{sec:heuristic}) to generate a set of feasible charging schedules. The best solutions are used to create initial CB cuts to restrict the MP's feasible region and we further strengthen the MP by adding cuts derived from the state-of-charge Constraints \eqref{eq:soc}. Following this initialization step, the main CB loop begins. In each iteration, we first solve the MP to get a candidate solution in terms of the binary $x$ and $y$ variables only. Assuming we obtain one, we solve the subproblem, which aims to find a solution that outperforms the current incumbent. If the SP is infeasible, indicating that this solution cannot outperform the incumbent, we use the CB cut generation procedure to exclude it. If it is feasible, our current MP solution becomes the new incumbent. We update it as such, then re-solve the SP and add CB cuts.

Sections \ref{sec:mp}-\ref{sec:cbcuts} next describe the individual steps of our Combinatorial Benders implementation in detail, including the MP and SP formulations as well as cut generation.

\subsection{Master Problem}
\label{sec:mp}

\begin{align}
    \min_{x^l_i, y^l_{ij}} \quad & \sum_{(l, i, j) \in \setA} c^l_{ij} y^l_{ij} & \label{eq:mpobj} \\
    \text{s.t.} \quad & \sum_{i \in \setT} y^l_{si} = 1 & \forall \enskip l \in \setC \label{eq:mpqueueconn1}\\
    & \sum_{i \in \setT} y^l_{it} = 1 & \forall \enskip l \in \setC \\
    & \sum_{j: (l, j, i) \in \setA} y^l_{ji} - \sum_{j: (l, i, j) \in \setA} y^l_{i j} = 0 & \forall \enskip l \in \setC, \enskip i \in \setT \label{eq:mpqueueconn3} \\
    & x^l_i = \sum_{j: (l, i, j) \in \setA} y^l_{ij} & \forall \enskip l \in \setC, \enskip i \in \setT \label{eq:mpxy} \\
    & \sum_{(i, j) \in \setA(\mathcal{S})} y^l_{ij} \leq \lvert \mathcal{S} \rvert - 1 & \forall \enskip \mathcal{S} \subset \setT: \lvert \mathcal{S} \rvert \geq 2 \label{eq:secs} \\
    & x^l_i \in \{0, 1\} & \forall \enskip l \in \setC, \enskip i \in \setT \\
    & y^l_{ij} \in \{0, 1\} & \forall \enskip (l, i, j) \in \mathcal{A} \label{eq:mpend}
\end{align}

The CB master problem \eqref{eq:mpobj}--\eqref{eq:mpend} includes all constraints on the binary variables $x$ and $y$ from the complete formulation, including the charger sequencing constraints \eqref{eq:mpqueueconn1}--\eqref{eq:mpqueueconn3} and $x$/$y$ relation \eqref{eq:mpxy}. The MP also introduces new constraints \eqref{eq:secs}, which ensure the MP solution does not contain any subtours. In Constraints \eqref{eq:secs}, $\setA(\mathcal{S})$ denotes the set of arcs contained in the trip subset $\mathcal{S}$, i.e.,

\begin{equation}
    \setA(\mathcal{S}) = \{(l, i, j) \in \setA: i \in S, j \in S\}
\end{equation}

Subtour elimination constraints were not necessary in the complete problem formulation \eqref{eq:mipobj2}--\eqref{eq:mipend} because any subtour would produce unbounded delay. But in the CB framework, that delay calculation is delegated to the subproblem, so the master problem will frequently generate solutions with subtours if Constraints \eqref{eq:secs} are omitted. We generate these constraints ``lazily'' within a branch-and-cut framework by checking for subtours and adding constraints only when subtours are detected, but having to handle these constraints directly in the MP is still a disadvantage of the CB method for our problem. 

Because the full problem objective function \eqref{eq:mipobj2} is not dependent on the binary variables, the MP essentially has a constant objective. In practice, this means that we may select any objective we wish, so we use a heuristic objective function to encourage the MP to generate integer solutions that are more likely to result in optimal delays as calculated by the subproblem. To do so, we define cost coefficients $c^l_{ij}$ for all arcs $(l, i, j) \in \setA$ as the lower bound on delay that would result from using that arc.

This lower bound is derived from the queuing and delay constraints. Let $j^+$ be the successor trip of $j$ on its same block, analogous to the previously defined $i^-$. If arc $(l, i, j)$ is used, then trip $j^+$ might be delayed due to bus $b_j$ queuing while $b_i$ charges. Specifically, $p_j \geq p_i$ by \eqref{eq:bigmfull} and because $t^l_i \geq 0$. Accounting for \eqref{eq:pluginlbfull}, we have $p_j \geq \sigma_i + \tau_i$, and then $d_{j^+} \geq \sigma_i + \tau_i - \sigma_{j^+}$. So, we set the MP cost for most arcs according to Equation \eqref{eq:mpcost}:

\begin{equation}
    c^l_{ij} = \max \left\{0, \sigma_i + \tau_i - \sigma_{j^+} \label{eq:mpcost} \right\}
\end{equation}

All dummy arcs (see Figure \ref{fig:fullnetworkfig}) and any arcs $(l, i, j)$ for which $j$ has no successor trip are assigned a cost of 0. Setting the costs with Equation \eqref{eq:mpcost} disincentivizes the solver from setting $y^l_{ij} = 1$ when doing so is guaranteed to delay trip $j^+$. In our experiments, this objective function gave better performance than a constant objective that essentially chose feasible solutions at random.

\subsubsection{Master Problem Cuts}
When initializing the MP, we also add some cuts to strengthen its formulation. Adding cuts helps to generate MP solutions that are more likely to be feasible for the SP, which decreases the total number of iterations that must be performed.

\textcolor{black}{We generate MP cuts based on the state-of-charge constraints \eqref{eq:socfull}, using them to derive lower bounds on the total number of times each bus must charge. Essentially, we convert Constraints \eqref{eq:socfull} to restrict the MP decision variables $x^l_i$ rather than the continuous variables $t^l_i$, which are now delegated to the subproblem.} To derive the MP cuts, we first note that whenever a bus visits any charger after trip $i$, the maximum amount of charge it can gain is $u^{\text{max}}_i$, the usable battery capacity of bus $b_i$. Now, for any trip $i$, let $l_i$ denote the charger located at the final stop of $i$, i.e., $l_i = l \in \setC: \overline{t^l_i} > 0$. Then we obtain valid inequalities \eqref{eq:mpcutstep1} based on Equations \eqref{eq:tmax} and \eqref{eq:chgtimefull}--\eqref{eq:socfull}:

\begin{equation}
    \sum_{j \in \theta(i)} u^{\text{max}}_i x^{l_i}_i \geq \alpha_i \quad \forall \enskip i \in \setT \label{eq:mpcutstep1}
\end{equation}

Dividing the inequalities \eqref{eq:mpcutstep1} by $u^{\text{max}}_i$ and adding the ceiling operator to the right-hand side since the left-hand side is always integral, we have:

\begin{equation}
    \sum_{j \in \theta(i)} x^{l_i}_i \geq \left\lceil \frac{\alpha_i}{u^{\text{max}}_i} \right\rceil \quad \forall \enskip i \in \setT \label{eq:mpcutfinal}
\end{equation}

\textcolor{black}{Equations \eqref{eq:mpcutfinal} yield one MP cut for every trip in the problem instance. In practice, we do not need to add all of these to the MP. Rather, we create these cuts by proceeding through each block in order of trip sequence and calculating the right-hand side $\left\lceil \frac{\alpha_i}{u^{\text{max}}_i} \right\rceil$. Each trip we progress to will add one variable $x^{l_i}_i$ to the left-hand side and may or may not increase the right-hand side, since $\alpha_i$ is nondecreasing as we progress through a block.  We therefore add a new cut to the MP each time the right-hand side increases to a new integer value; if it does not increase, then the cut for trip $i$ is dominated by the one for $i^-$ and does not improve the MP formulation.}

\subsection{Subproblem}
Solving the CB master problem provides a candidate solution in terms of the $x$ and $y$ variable values. The subproblem's role is then to check whether it could be an optimal solution to the complete problem by verifying both that it is a feasible choice of integer variable values and that it can minimize delay. This optimality criterion is evaluated by tracking an incumbent solution and its total delay value, denoted $z^*$, as described in \citet{Codato2006}.

For a given MP iteration $k$, let $\setA_1^{(k)}$ be the set of arcs used in the solution, i.e., $(l, i, j) \in \setA_1^{(k)}$ if $y^{l(k)}_{ij}=1$. Likewise, let $\setT^{(k)}_0$ be the set of trip-charger pairs for which charging is skipped, i.e., $(l, i) \in \setT^{(k)}$ if $x_i^{l(k)} = 0$. Then the SP for iteration $k$ is:

\begin{align}
    \min_{d_i, p_i, t^l_i} \quad & \sum_{i \in \setT} d_i & \label{eq:spobj} \\
    \text{s.t.} \quad & \sum_{i \in \setT} d_i \leq z^* - \epsilon & \label{eq:spincumbent} \\
    & d_{i} \geq p_{i^-} + t^l_{i^-} - \sigma_{i} & \forall \enskip l \in \setC, \enskip i \in \setT: \enskip n_i >= 2 \\
    & p_{i} \geq \sigma_{i} + d_{i} + \tau_{i} & \forall \enskip i \in \setT \\
    & p_{j} \geq p_{i} + t^l_{i} & \forall \enskip (l, i, j) \in \setA_1^{(k)} \label{eq:sppluginqueue} \\
    & \alpha_i \leq \sum_{j \in \theta(i)} \sum_{l \in \setC} \rho^l t^l_j \leq \beta_i & \forall \enskip i \in \setT \\
    & t^l_{i} = 0 & \forall \enskip (l, i) \in \setT^{(k)}_0 \label{eq:spzerocharge} \\
    & 0 \leq t^l_{i} \leq \overline{t^l_{i}} x^l_i & \forall \enskip l \in \setC, \enskip i \in \setT \\
    & d_i \geq 0 & \forall \enskip i \in \setT \label{eq:spend}
\end{align}

The CB subproblem \eqref{eq:spobj}--\eqref{eq:spend} aims to minimize the delay based on the binary decisions from the MP. Constraint \eqref{eq:spincumbent} requires that the subproblem solution must improve the current incumbent objective value $z^*$ by some small tolerance $\epsilon$, which ensures that when the MP is found to be infeasible, the current incumbent is optimal. Constraints \eqref{eq:sppluginqueue} and \eqref{eq:spzerocharge} replace the corresponding big-M constraints from the full problem, enforcing bounds on plugin time $p_i$ and charging time $t^l_i$ only as needed based on the current MP solution. The remaining SP constraints are identical to those from the complete formulation in Section \ref{sec:formulation}.

\subsection{Combinatorial Benders Cuts}
\label{sec:cbcuts}
As shown in Figure \ref{fig:cbflowchart}, when the SP is infeasible, we add Combinatorial Benders cuts to exclude the incumbent MP solution. The key to creating these cuts is to identify a minimal infeasible subsystem (MIS), also called an irreducible infeasible subsystem (IIS) of the infeasible SP instance. An MIS is an inclusion-minimal set of rows of the SP constraint matrix, where ``inclusion-minimal'' means that if any single constraint is removed, the resulting subsystem of constraints admits a feasible solution. \citet{Codato2006} describe an algorithm to generate multiple such MISs for any MP solution that yields an infeasible SP.

An MIS of the subproblem identifies a set of MP variables---corresponding to the instances of the conditional constraints \eqref{eq:sppluginqueue} and \eqref{eq:spzerocharge} included in the MIS---that have forced the problem to be infeasible. To make sure that the MP does not produce a solution with this same set of variables again, we add a CB cut to exclude it. For this purpose, we can represent the MIS as two sets of indices: one set $\setM_x$ corresponding to Constraints \eqref{eq:sppluginqueue} and another $\setM_y$ corresponding to Constraints \eqref{eq:spzerocharge}. That is, $\setM_x = \left\{(l, i) \in \setM: x^l_i = 0 \right\}$ and $\setM_y = \left\{(l, i, j) \in \setM: y^l_{ij} = 1 \right\}$. Given such an MIS, a CB cut is formed by the simple inequality \eqref{eq:cbcut}:

\begin{equation}
    \sum_{(l, i) \in \setM_x} x^l_i + \sum_{(l, i, j) \in M_y} (1 - y^l_{ij}) \geq 0 \label{eq:cbcut}
\end{equation}

Equation \eqref{eq:cbcut} enforces that at least one of the binary variables $x^l_i$ and $y^l_{ij}$ included in the MIS must change its value in order to obtain a feasible solution.

\subsection{Implementation Details}
We implemented the CB algorithm using Python and the Gurobi solver via the \texttt{gurobipy} package. Following Figure \ref{fig:cbflowchart}, we begin by initializing the master problem and adding CB cuts for heuristic solutions. In our experiments, we generate cuts for any heuristic solution with an objective value within 50\% of the best identified objective. Once the MP is initialized, we run the CB algorithm within a branch-and-cut framework using Gurobi's callback capabilities. When a new optimal solution to the MP is detected, we first check if it contains any subtours. If it does, we cut off any such subtours by adding Constraints \eqref{eq:secs} as lazy constraints. If it does not, we progress to solving the subproblem and generating CB cuts. To generate an MIS each time that the SP is infeasible, we use Gurobi's \texttt{computeIIS()} function, which identifies a single MIS out of many possibilities. Identifying multiple MISs (and consequently multiple cuts) can help the algorithm converge faster, so each time an MIS is found, we remove one constraint from that MIS at random, verify that the relaxed SP model is still infeasible, and run \texttt{computeIIS()} again to find a new unique MIS. We repeat this procedure until the relaxed SP model is feasible, which typically results in finding several MISs per CB iteration.

\section{Heuristic Solution Method: Select--Sequence--Schedule}
\label{sec:heuristic}
Deriving the CB algorithm for our charge scheduling formulation \eqref{eq:mipobj2}--\eqref{eq:mipend} inspired a heuristic algorithm that follows a similar pattern, but with much faster convergence. The heuristic design is based on a few key insights presented by the decomposition approach. First, most of the problem's complexity comes from the two interconnected decisions of selecting which trips include charging (the values of $x^l_i$) and the sequence in which these trips are served by each charger (the values of $y^l_{ij}$). Once these values are set, it is easy to determine the optimal charging durations and complete schedule of plugin/departure/delay times with the subproblem \eqref{eq:spobj}--\eqref{eq:spend}, which is just an LP. Additionally, the problem of selecting the order in which to serve trips may not be too difficult in practice once the charging trips have been selected. We can expect that the optimal charging order is unlikely to differ too much from simple first-in, first-out priority.

Based on this logic, we devised a heuristic algorithm based on relaxing and separating the MIP formulation \eqref{eq:mipobj2}--\eqref{eq:mipend}. We relax the complicating queue constraints \eqref{eq:bigmfull} and restructure the problem into three phases we call \textit{selection}, \textit{sequencing}, and \textit{scheduling}, together forming what we call the 3S algorithm.

\begin{figure}[ht]
    \begin{center}
    \begin{tikzpicture}[node distance=2cm]
    \node (getcosts) [start] {Generate Costs $\theta^l_i$};
    \node (phase1) [process, below of=getcosts, text width=5cm] {\textbf{Phase 1}\\\textit{Select} trips and chargers};
    \node (phase2) [process, below of=phase1, text width=5cm] {\textbf{Phase 2}\\\textit{Sequence} trips served by each charger};
    \node (phase3) [process, below of=phase2, text width=5cm] {\textbf{Phase 3}\\ \textit{Schedule} events precisely};
    \node (outputs) [stop, below of=phase3] {Charge durations $t^l_i$, plugin times $p_i$, delays $d_i$};
    
    \draw [arrow] (getcosts) -- (phase1);
    \draw [arrow] (phase1) -- (phase2);
    \draw [arrow] (phase2) -- (phase3);
    \draw [arrow] (phase3) -- (outputs);
    \end{tikzpicture}
    \end{center}
    \caption{Flowchart of 3S heuristic solution procedure.}
    \label{fig:heurflowchart}
\end{figure}
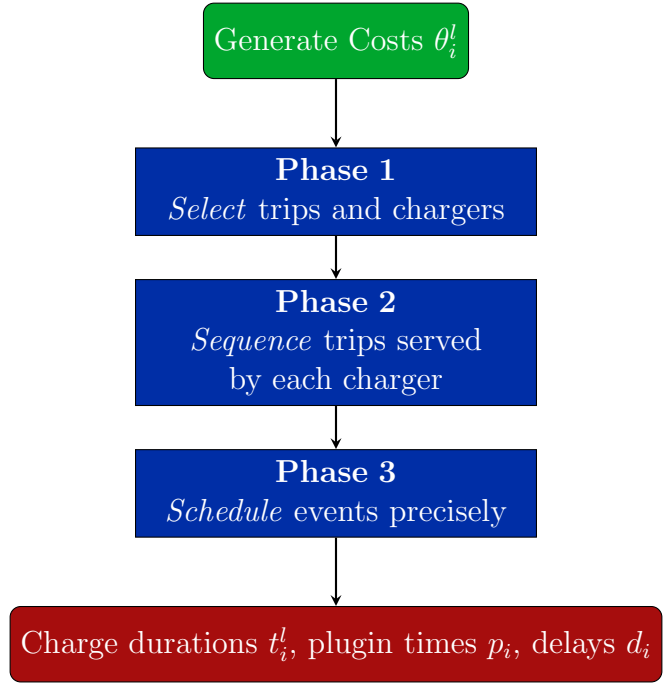

Figure \ref{fig:heurflowchart} gives an overview of the 3S heuristic algorithm. The first step is to initialize some random cost parameters $\theta^l_i$, which helps the algorithm explore a wider range of feasible solutions. Then, in the selection phase, we solve a separate LP for each bus to select a set of trips when it will charge. This is equivalent to setting the values of $x$ or populating the set $\setT^{(k)}_{0}$ in the CB approach. In the sequencing phase, we perform a simple sorting operation for each charger to set the order in which it serves trips. This is analogous to setting the $y$ variable values or the set $\setA^{(k)}_1$ from the CB master problem. In the final \textit{scheduling} phase, we solve one more LP to optimize charging durations and the resulting delay given the selection and sequence decisions. This final LP is identical to the CB subproblem with the incumbent bound constraint \eqref{eq:spincumbent} removed. The final output of any run of the 3S heuristic is a feasible solution to \eqref{eq:mipobj2}--\eqref{eq:mipend} and its corresponding objective value. Sections \ref{sec:selection}--\ref{sec:scheduling} next describe each of the three phases in detail.

\subsection{Phase 1: Trip Selection}
\label{sec:selection}
The purpose of Phase 1 is to set feasible values of the binary charging decisions $x^l_i$ for all trips $i$ and chargers $l$. To do this, we relax Constraints \eqref{eq:bigmfull}, which are used to accurately calculate the delays $d_i$ but do not affect feasibility with respect to state of charge. When Constraints \eqref{eq:bigmfull} are removed from the MIP formulation \eqref{eq:mipobj2}--\eqref{eq:mipend}, the binary variables $y$ are no longer needed since they have no impact on the objective function. The binary variables $x$ are also no longer needed, since these were just dummy variables used to link $y$ to the charging durations $t$. Finally, we can absorb the plugin time variables $p_i$ into $d_i$ now that the plugin time is controlled only by Constraints \eqref{eq:pluginlbfull}. With these variables and the constraints that involve them removed from the full formulation, we arrive at the Phase 1 LP \eqref{eq:p1obj}--\eqref{eq:p1end}:

\begin{align}
    \min_{d_i, t^l_i} \quad & \sum_{i \in \setT} \left[d_i + \sum_{l \in \setC} \theta^l_i t^l_i \right] & \label{eq:p1obj} \\
    \text{s.t.} \quad & d_{i} \geq \sigma_{i^-} + d_{i^-} + \tau_{i^-} + t_{i^-} - \sigma_{i} & \forall \enskip i \in \setT: \enskip n_i \geq 2 \\
    & \alpha_i \leq \sum_{j \in \theta(i)} \sum_{l \in \setC} \rho^l t^l_j \leq \beta_i & \forall \enskip i \in \setT \\
    & 0 \leq t^l_{i} \leq \overline{t^l_{i}} & \forall \enskip i \in \setT, \enskip l \in \setC \label{eq:p1battery} \\
    & d_i \geq 0 & \forall \enskip i \in \setT \label{eq:p1end}
\end{align}

The Phase 1 objective function \eqref{eq:p1obj} includes a randomization term $\sum_{i \in \setT, l \in \setC} \theta^l_i t^l_i$, where $\theta^l_i$ is randomly sampled from a suitable distribution. This term helps to generate a variety of feasible solutions through repeated runs of the 3S algorithm. Since there will often be many different charging patterns that result in the same optimal delay in Phase 1 (especially if the minimum total delay is zero), adding this random coefficient to the charging times helps 3S to consider a greater variety of potentially good first-stage solutions. A distribution should be chosen so that $\theta^l_i$ takes on a variety of both positive and negative values. Negative values of $\theta^l_i$ are useful to produce results in which a bus distributes its charging over a greater number of trips, which may help to reduce queuing. It is also important that the values of $\theta^l_i$ are not so large that the randomization term dominates the delay term in the objective function, which could result in poor performance. In our experiments, we sample the values of $\theta^l_i$ from a normal distribution with mean 0 and standard deviation 0.5.

The Phase 1 LP \eqref{eq:p1obj}--\eqref{eq:p1end} is separable across all buses in the system, since the conditional plugin time constraints \eqref{eq:bigmfull} that linked trips on different blocks have been removed. Phase 1 therefore consists of $N_b$ independent linear programs, where $N_b$ is the number of buses serving the trips $\setT$. Note that since the number of variables only scales with the number of trips served by a single vehicle, each LP is quite small.

Solving the Phase 1 LP \eqref{eq:p1obj}--\eqref{eq:p1end} yields first-stage charging durations $t^l_i$ and delays $d_i$ for all trips $i$ and chargers $l$. We convert the charging durations to their equivalent binary values $x^l_i$, i.e. $x^l_i=1$ if $t^l_i > 0$ and $x^l_i=0$ otherwise. These binary decisions also correspond to the set of non-charging trips $\setT^{(k)}_0$ from the CB master problem. Each trip $i$ for which $\sum_{l \in \setC} x^l_i=1$ is also assigned a first-stage plugin time $p_i=d_i + \sigma_i + \tau_i$ when it is expected to start charging. These form the inputs to Phase 2.

\subsection{Phase 2: Charger Sequencing}
\label{sec:sequencing}
In Phase 2, we process the charging decisions from Phase 1 into a sequence for every charger that is used. For each charger $l$, we simply sort all trips $i$ for which $x^l_i=1$ by their first-stage plugin times $p_i$. Sorting these trips gives us a simple mapping to a feasible set of charging arcs $\setA_1^{(k)}$: for any consecutive pair of trips $i$ and $j$, we add arc $(l, i, j)$ to $\setA_1^{(k)}$. For each charger we also add arcs connecting to the initial and final dummy nodes, respectively. That is, if $i$ is the first and $j$ is the last trip served by charger $l$, we add the arcs $(l, s, i)$ and $(l, j, t)$.

\subsection{Phase 3: Event Scheduling}
\label{sec:scheduling}
Phase 3 takes the results obtained from Phases 1 and 2---essentially, all the binary decisions of the complete MIP formulation---and uses them to define an exact charging schedule to minimize delays. Since the first two phases populated the sets $\setT_0^{(k)}$ and $\setA_1^{(k)}$, we have all the inputs necessary to construct an instance of the subproblem \eqref{eq:spobj}--\eqref{eq:spend}, though we need to either set the incumbent delay value $z^*$ to an arbitrarily high value or remove Constraint \eqref{eq:spincumbent} entirely. Since the subproblem is a linear program whose number of variables and constraints scales roughly with the number of trips in an instance, it can generally be solved quite quickly. Solving the subproblem yields final values of charging duration and delay for each trip that account for queuing. Note that in Phase 3 we allow for changing the charging durations originally output by Phase 1; the purpose of Phase 1 is to make the binary decisions of when and where to charge, but we relax the problem to an LP for efficiency. Also note that including the charging time constraints \eqref{eq:p1battery} in the Phase 1 LP for each bus ensures that Phase 3 is always feasible; the purpose of Phase 3 is to optimize charging durations to exactly determine the minimum possible delay for the given set of charging trips and sequence for each charger.

\section{Case Studies}
\label{sec:casestudy}
We constructed a variety of instances in order to analyze the performance of the CB method and 3S heuristic. These instances were based on two different transportation networks. The first, described in Section \ref{sec:simplenetwork}, is a small notional network consisting of two bus routes served by a single shared charger. The small network allows us to evaluate the performance (in terms of solution time and optimality gap) of both our methods in comparison to directly solving the model with Gurobi. The second network, described in Section \ref{sec:metronetwork}, is based on the actual transit system operated by King County Metro in the greater Seattle area. This example shows how the 3S heuristic can support transit operations at a real-world scale.

\subsection{Simple Notional Network}
\label{sec:simplenetwork}
\begin{figure}[ht]
    \centering
    \begin{tikzpicture}
        \Vertex[label=\textcolor{tlgreen}{\large \Lightning}, x=0,y=0, size=0.5, color=white]{s}
        \Vertex[x=2,y=1, size=0.5, color=white]{e}
        \Vertex[x=0, y=2, size=0.5, color=white]{n}

        \node[orange, rotate=90] at (-0.7, 1) {\large Route A};
        \Edge[Direct, bend=15, color=tlblue](s)(e)
        \Edge[Direct, bend=15, color=tlblue](e)(s)
        \Edge[Direct, bend=15, color=orange](s)(n)
        \Edge[Direct, bend=15, color=orange](n)(s)
        \node[tlblue, rotate=25] at (1.4, 0) {\large Route C};

    \end{tikzpicture}
    \caption{Simple case study network.}
    \label{fig:simplecase}
\end{figure}
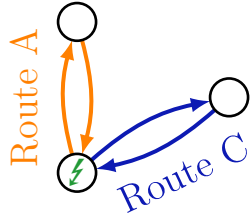

To test the CB method and assess the performance of the 3S heuristic, we used a small test network originally presented in Appendix B of \citet{mccabe_optimal_2023}. That original network was further simplified to include only two of the routes (A and C) and a single charger at their shared terminal, as sketched in Figure \ref{fig:simplecase}. Route A's headway was decreased from 20 to 30 minutes to limit the number of vehicles in the case study, but all other parameter values remained the same as reported in \citet{mccabe_optimal_2023}. All buses have 400 kWh of usable battery capacity and consume energy at a rate of 3 kWh/mi (1.86 kWh/km); we test a variety of charger power levels from 300 to 500 kW. Table \ref{tab:simplerouteparams} documents the schedule and distance parameters of each route. With such a schedule, the network consists of 8 buses that complete 84 total trips. 

\begin{table}[ht]
    \centering
    \caption{Basic parameters for notional routes used in simple case study.}
    \begin{tabular}{l|c|c}
        Route & A & C \\ \hline
        Distance (mi.) & 15 & 15 \\
        Time (min.) & 40 & 45 \\
        Layover Time (min.) & 20 & 15 \\
        Headway (min.) & 30 & 60 \\
        One-Way Distance & 15 mi. (24 km) & 15 mi. (24 km)
    \end{tabular}
    \label{tab:simplerouteparams}
\end{table}

\subsection{King County Metro Network}
\label{sec:metronetwork}
The second test network is based on the transit system of South King County, WA, USA. We use this network to study the performance of the 3S heuristic, as it is too large to be solved to optimality by either Gurobi or our CB method. The case study includes some of the busiest routes planned for electrification in the near future: the RapidRide lines F and H as well as routes 131, 132, 150, 153, and 161. We collected relevant data including trip schedules (corresponding to $\sigma_i$ and $\tau_i$ values) and distances as well as block sequences from the Metro GTFS feed \citep{MetroGTFS} that was posted November 27, 2023. We used Wednesday, December 6, 2023 as a test date.

\subsubsection{Data Collection}
To filter down the GTFS data to the scope of our case study, we first identified all \texttt{service\_id} values active on the case study date and all blocks and trips active for these service IDs. To set the cumulative charging bounds $\alpha_i$ and $\beta_i$ for each trip $i$, we first determined the distance of each trip by matching it to its corresponding \texttt{shape\_id} and the maximum \texttt{shape\_dist\_traveled}. We assumed each bus had a usable battery capacity of 300 kWh and an average energy consumption rate of 3 kWh/mi (1.86 kWh/km) for all trips. The first and last trips of each block also had increased energy demand due to pull-out and pull-in trips from/to the depot. We assumed all buses were kept overnight at Metro's South Base and calculated pull-out/pull-in distances based on driving directions between the base and all relevant terminal locations, as calculated by the Openrouteservice API \citep{Openrouteservice}. To set the upper bound on charging time $\overline{t^l_i}$ for all trips and chargers, we identified the final stop coordinates of each trip and calculated the driving distance to each charger (again using Openrouteservice) and assumed buses were able to use a charger if the distance to it was less than 0.25 miles (0.4 km). 

\begin{figure}[ht]
    \centering
    \includegraphics[width=0.45\linewidth]{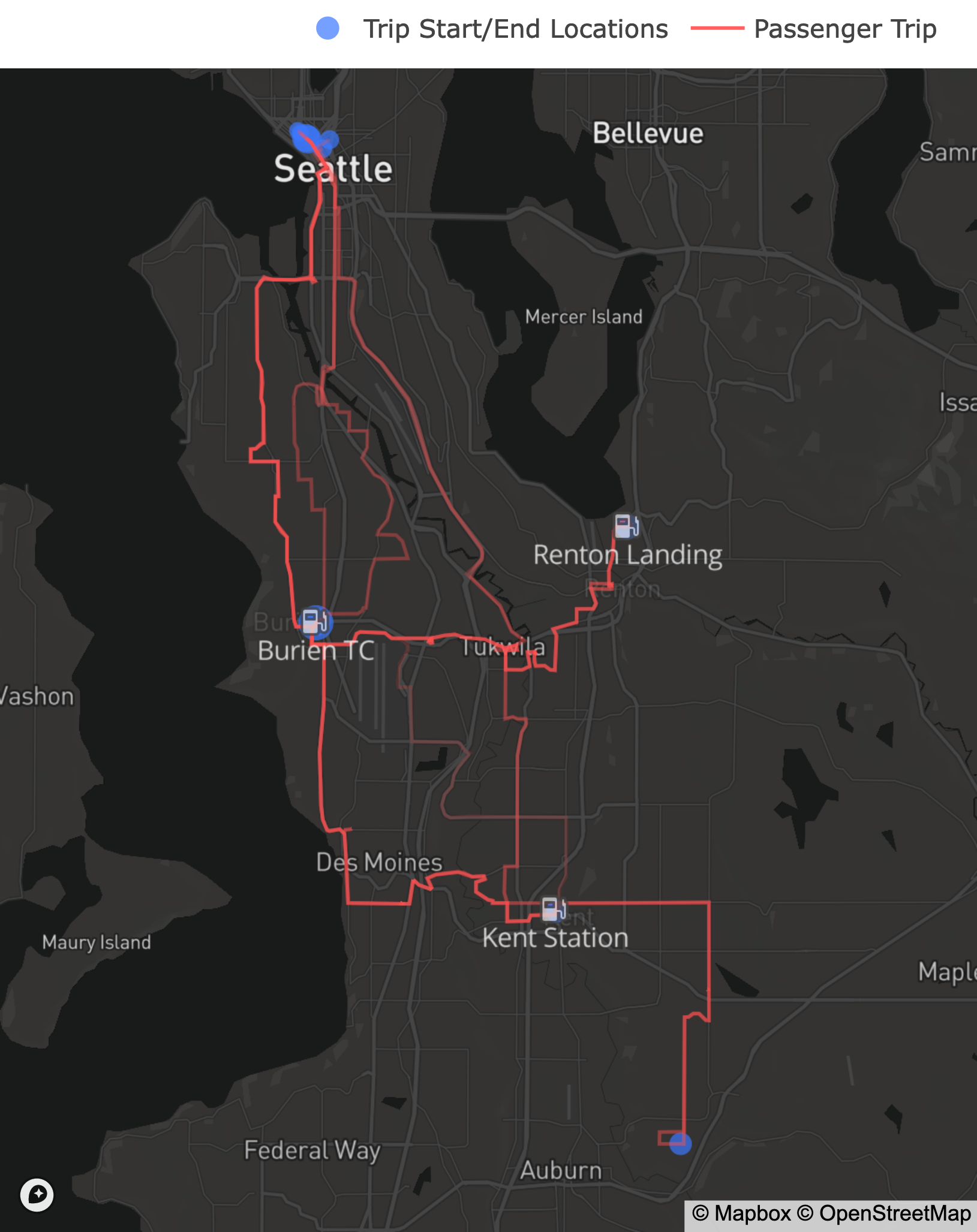}
    \caption{Map showing trips, terminals, and charger locations for the King County case.}
    \label{fig:instancemapb}
\end{figure}

The case study includes all blocks that serve trips on the RapidRide F Line and H Line as well as routes 131, 132, 150, 153, 161, and 165. Buses on these routes are served by one charger each at the Burien Transit Center (which serves the F Line, H Line, 131, 132, and 165), Renton Landing (the eastern terminal of the F Line), and Kent Station (the southern terminal of routes 150, 153, and 161). We assumed all chargers had an identical power output of 450 kW in each test case. Figure \ref{fig:instancemapb} maps all trips included in this scope as well as the three chargers considered.

In total, the King County network we analyze includes 42 buses that need to use fast chargers at some point during the day. An additional 17 blocks that serve these routes are excluded because they can be completed with depot charging alone. The 42 fast-charging buses complete a total of 493 trips on the test date.

\subsubsection{Exogenous Delay Scenarios}
We introduce two scenarios to model exogenous delays, which highlight the flexibility of our modeling approach and the 3S algorithm. In Scenario A, all trips run exactly on schedule. In Scenario B, the trip duration $\tau_i$ is increased by 40\% for all trips with a scheduled departure between 7:00 a.m. and 9:00 a.m., which is meant to roughly simulate heavy traffic during these commute hours. In this case, buses are delayed enough that charging within scheduled layover time is impossible, and we search for a solution that yields the minimum possible total delay. The results in Section \ref{sec:results} show how our model adapts by shifting charging later in the day, once buses have resumed operating on schedule.

\subsection{Results}
\label{sec:results}
\subsubsection{Simple Network}
We solved the model for the simple network of Section \ref{sec:simplenetwork} at a variety of different charger power levels from 300 to 500 kW. The two exact methods were subject to a time limit of one hour for each instance; the 3S heuristic was run 500 times for each case, though it was terminated as soon as a solution with zero delay was found, if applicable. The results are summarized in Table \ref{tab:simpleresults}. For each method, we report the objective value of the best solution obtained, the time to find the first solution with that best objective value, and the total solution time. A dash is used to represent that the algorithm timed out.

\begin{table}[ht]
    \centering
    \caption{Summary of results on simple case study instances. \textbf{BO}: best objective value (minutes of total delay). \textbf{T-BO}: time to find solution with best objective (s), \textbf{T-T}: total solution time (s). A dash indicates that the algorithm did not terminate.}
    \begin{tabular}{c|c|c|c|c|c|c|c|c|c}
        & \multicolumn{3}{c|}{\textit{\textbf{Direct Solve}}}
        & \multicolumn{3}{c|}{\textit{\textbf{3S Heuristic}}} & \multicolumn{3}{c}{\textit{\textbf{Benders}}} \\
        \textbf{Power} & \textbf{BO} & \textbf{T-BO} & \textbf{T-T} & \textbf{BO} & \textbf{T-BO} & \textbf{T-T} & \textbf{BO} & \textbf{T-BO} & \textbf{T-T}\\ \hline
        300 kW & 704 & 133 & - & 759 & 0.7 & 12 & 759 & N/A & -\\
        350 kW & 71.4 & 15 & 57 & 73.7 & 0.05 & 9 & 71.4 & 1463 & - \\
        400 kW & 5 & 56 & 66 & 5 & 0.7 & 10 & 5 & 0.7 & 47 \\
        450 kW & 0 & 55 & 55 & 0 & 0.1 & 0.1 & N/A & N/A & N/A \\
        500 kW & 0 & 138 & 138 & 0 & 0.4 & 0.4 & N/A & N/A & N/A
    \end{tabular}
    \label{tab:simpleresults}
\end{table}

In Table \ref{tab:simpleresults}, we see that optimal delay decreases and eventually reaches zero as the power level increases. For these small instances involving only 8 buses, Gurobi solved the problem to optimality in 4 out of 5 cases. The comparative performance of the heuristic and CB methods varied depending on the instance. At the two lower power levels which represent more difficult instances, the Benders algorithm failed to prove optimality within the time limit, though in the 350 kW case it did successfully identify an optimal solution with 71.4 minutes of delay. With 400 kW chargers, the CB method outperforms the direct solution via Gurobi. As intended, the 3S algorithm quickly identifies an optimal solution with 5 minutes of total delay and the CB procedure proves its optimality in 47 seconds, whereas directly solving the problem took 66 seconds. For the 450 and 500 kW cases, the direct solution approach struggles to find a good solution, whereas 3S identifies a zero-delay solution very quickly and the entire CB procedure can be skipped because total delay can never be negative.

\begin{figure}
    \centering
    \includegraphics[width=0.9\linewidth]{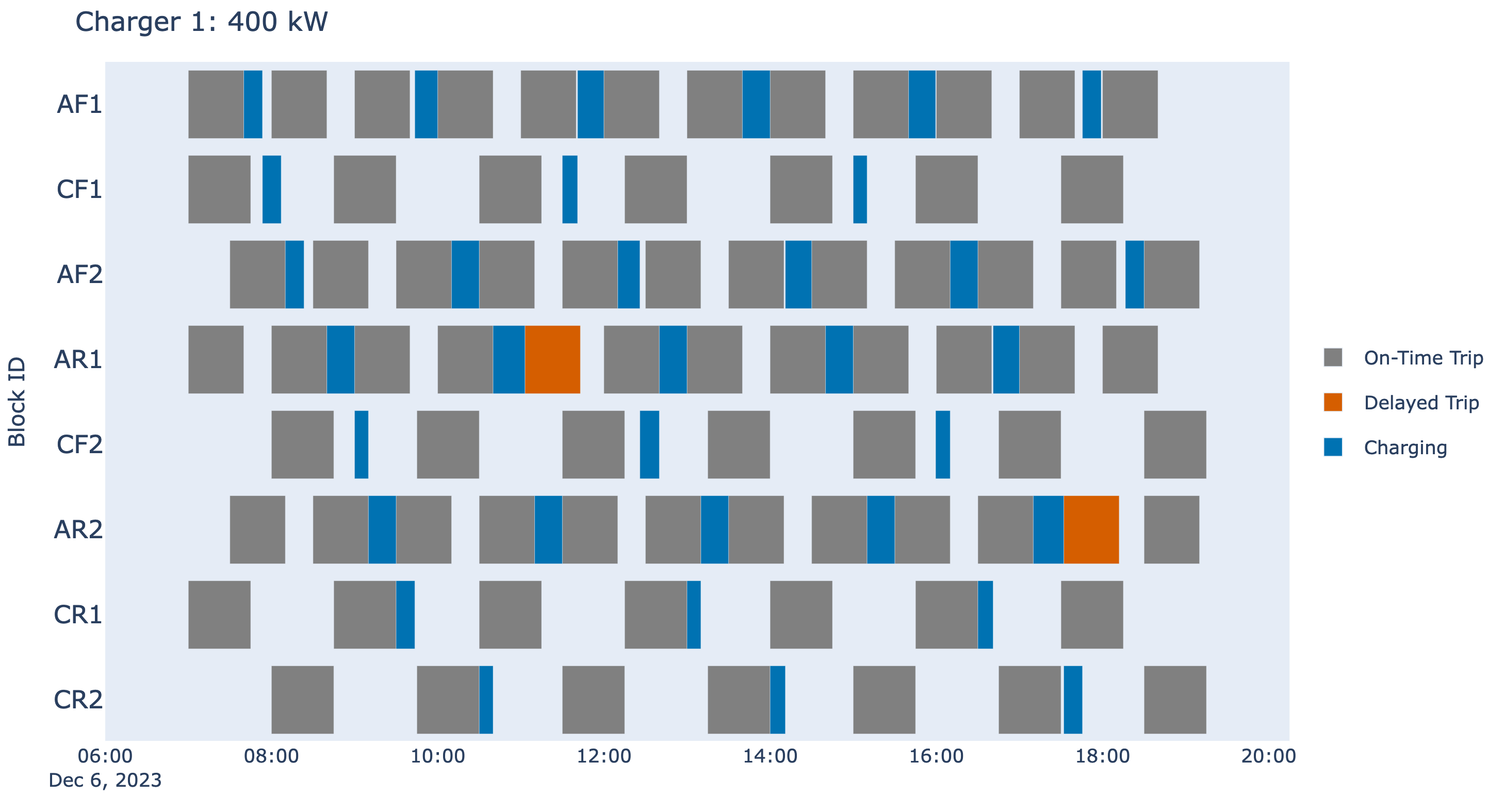}
    \caption{Optimized timeline of operations for the simple network with 400 kW chargers.}
    \label{fig:simpletimeline}
\end{figure}

To illustrate the results in more detail, Figure \ref{fig:simpletimeline} shows the complete timeline of activities for all buses in the 400 kW instance. Gray blocks on the timeline indicate that a bus is completing a passenger service trip that departed on time; delayed trips are shown in orange. A blue block in a bus's timeline indicates that it is plugged in at the charger. Looking at block CF1 for example, which completes 7 trips on Route C, we can see that it chargers after trips 1, 3, and 5, and all trips are on time. Blocks AR1 and AR2, with both complete more trips and require more charging, each have one trip that departs a few minutes late. We can see that this optimal solution involves frequent short charges---in fact, all eight buses buses plug in every time they are able to (every other trip).

These results highlight both the benefits and limitations of our methods. First, we note that the 3S heuristic identifies a near-optimal (often optimal) solution very quickly in all cases. Its worst performance was for the 300 kW case, when the best solution it identified had only 7.8\% more delay than the best one found by Gurobi. One key advantage of the 3S method is that its solution time does not depend significantly on the difficulty of the instance being solved; while neither of the exact methods converged within the time limit for the 300 kW instance, the heuristic produced good feasible solutions just as quickly as it did for other instances.

CB produced mixed results on our instances. We see that for the 400 kW instance, the CB method converges more quickly than solving directly, but it is much slower than Gurobi for the two harder instances. In the 350 kW case, it is notable that CB is able to improve the incumbent solution supplied by the heuristic, but it still does not terminate before the time limit. With a 450 or 500 kW charger, integrating Benders with the 3S approach means that we can avoid running the Benders process entirely since 0 is the best possible feasible delay.

There are a few likely reasons for the limited benefits of the CB approach for our charge scheduling model. First, delegating all scheduling and delay calculation to the subproblem limits the information that is available to the master problem, which limits Gurobi's ability to speed up the solution process via presolve methods and cuts. Based on our results, it seems that the cuts generated by the CB approach are not particularly strong compared to the cuts Gurobi generates on its own when provided with the complete formulation. Second, the CB method is primarily used for problems in which the objective depends only on the integer variables. While we followed the recommendations from \citet{Codato2006} on adapting the method to a problem where the objective depends instead on continuous variables, we found that the approach is not as well suited to this type of problem. For one, convergence is inevitably slower compared to the typical CB application because it usually requires a large number of CB cuts to determine that the MP is infeasible and terminate the algorithm. Additionally, because the MP objective has no real meaning for the complete problem, we do not get any information on the optimality gap. Whereas in the typical CB application the Benders cuts are embedded within a branch-and-cut procedure and the search process continually improves upper and lower bounds on the optimal objective value, in our case we have an upper bound from the incumbent solution, but do not obtain any lower bound until the MP is infeasible and optimality is proven.

\subsubsection{King County Metro Network}

\begin{table}[ht]
    \centering
    \caption{Summary of results on King County Metro network. \textbf{BO}: best objective value (minutes of total delay). \textbf{T-BO}: time to find solution with best objective (s), \textbf{T-T}: total solution time (s), \textbf{ND}: number of delayed trips.}
    \begin{tabular}{c|c|c|c|c}
        \textbf{Instance} & \textbf{BO} & \textbf{T-BO} & \textbf{T-T} & \textbf{ND} \\ \hline
        A & 0 & 0.32 & 0.32 & 0 \\
        B & 186.75 & 0.47 & 71 & 36 \\
    \end{tabular}
    \label{tab:metroresults}
\end{table}

Table \ref{tab:metroresults} documents the key results from applying the 3S heuristic to the King County network for the two delay scenarios. In Scenario A, the algorithm identifies an optimal solution with zero delay in a fraction of a second. In Scenario B where some delay is inevitable, it takes only 0.47 seconds to identify the best solution and just over a minute to run 500 iterations of 3S. 36 out of the 446 total trips are delayed in the best solution, with a maximum delay of 16.5 minutes for a single trip.

\begin{figure}
    \centering
    \includegraphics[width=0.9\linewidth]{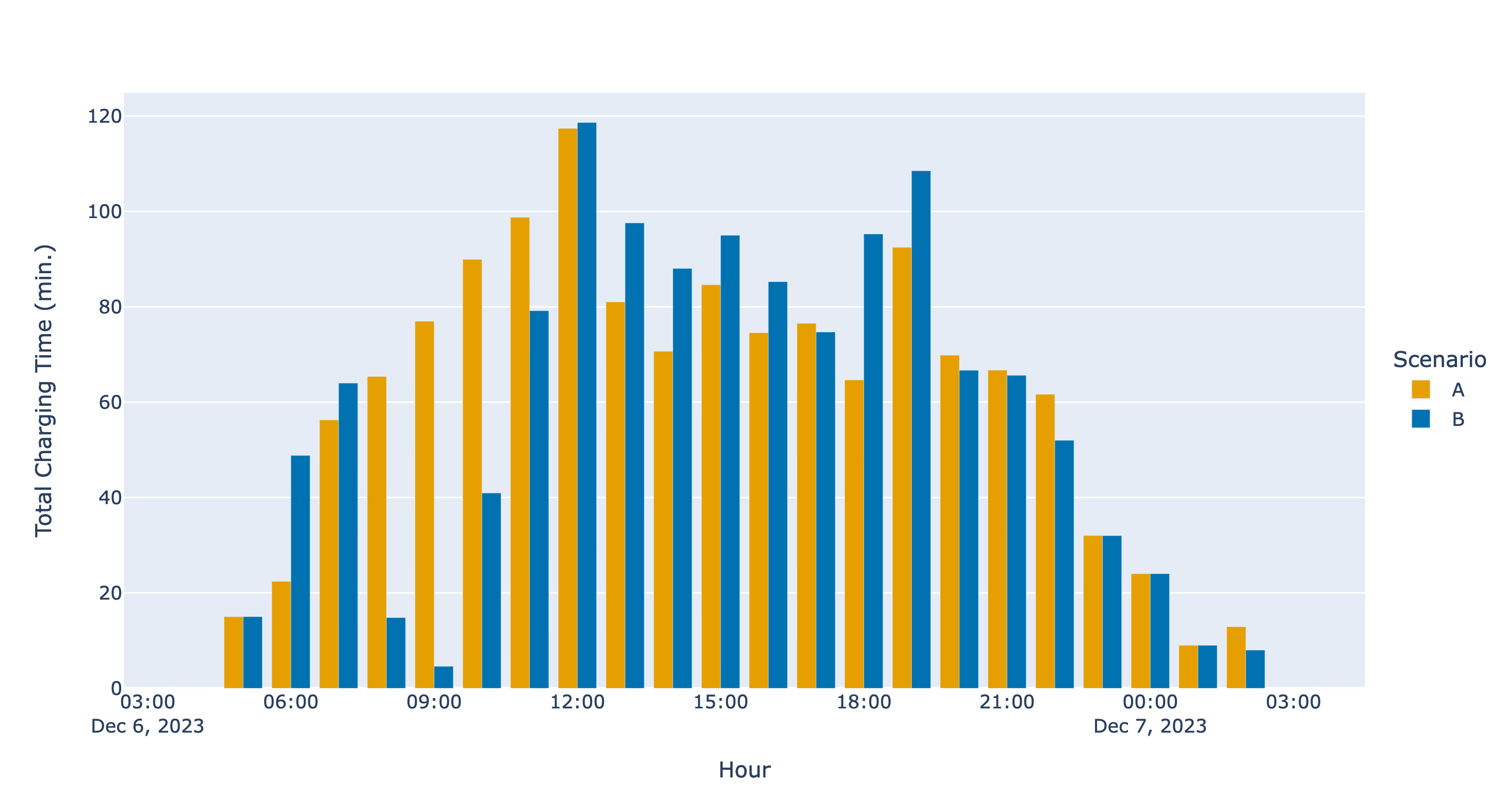}
    \caption{Total charging time by hour across all sites for the two scenarios.}
    \label{fig:hourlycomparison}
\end{figure}

\begin{figure}
    \centering
    \begin{subfigure}{\linewidth}
        \centering
        \includegraphics[width=\linewidth]{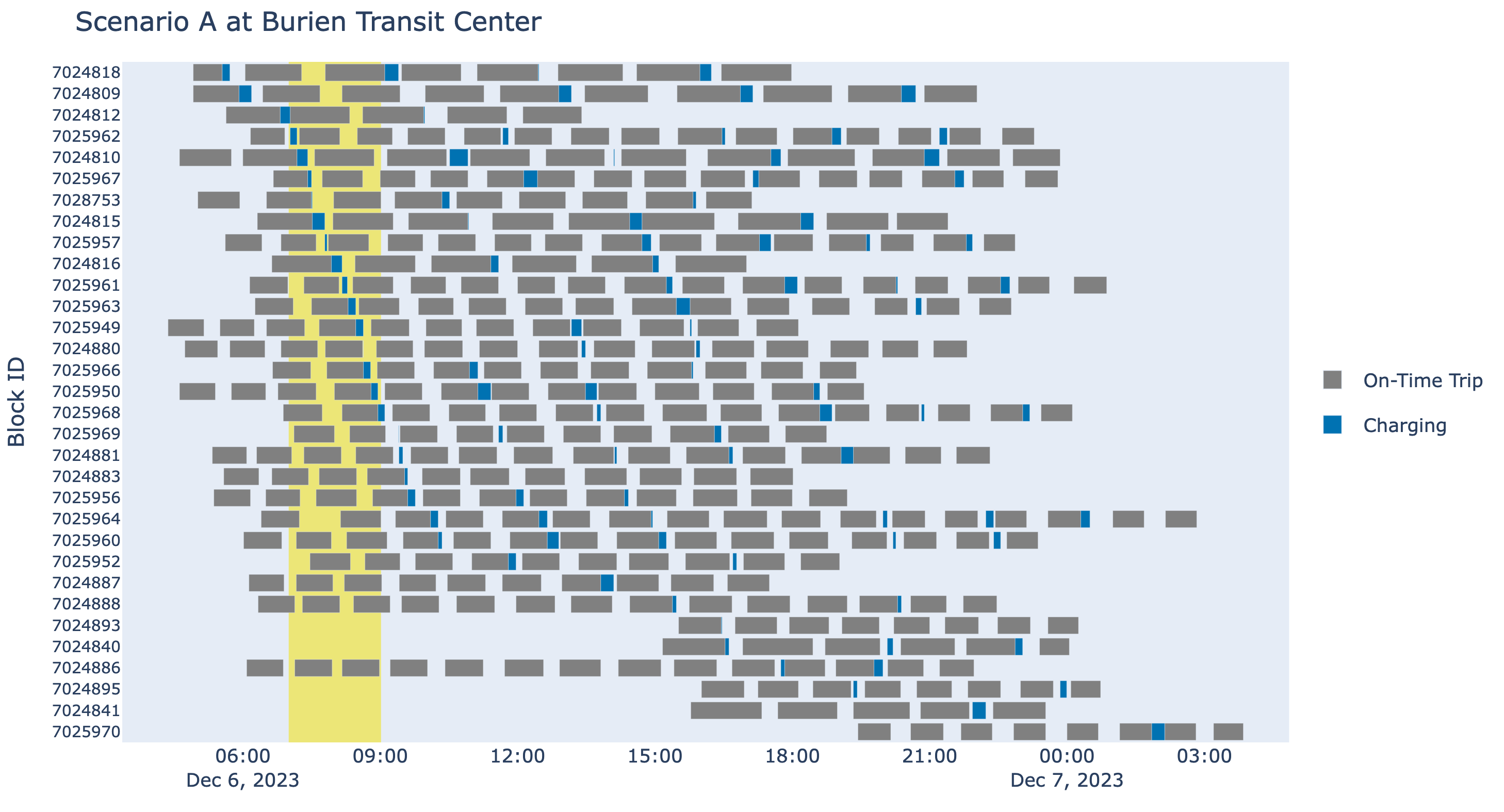}
    \end{subfigure}

    \begin{subfigure}{\linewidth}
        \centering
        \includegraphics[width=\linewidth]{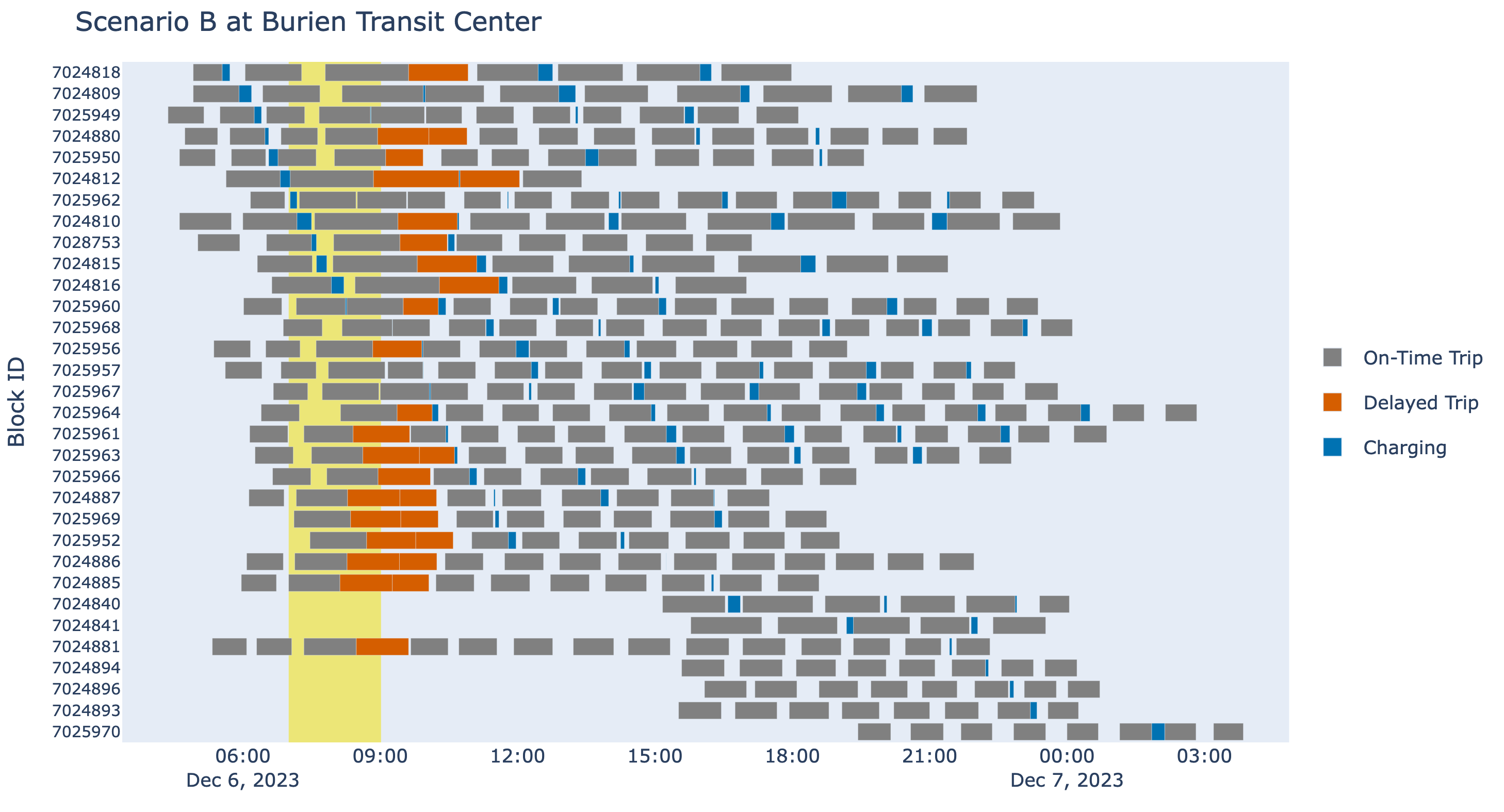}
    \end{subfigure}
    \caption{Timeline of scheduled charger utilization for Scenario A.}
    \label{fig:burientimeline}
\end{figure}

To show how the 3S algorithm adapts to our simulated delays in Scenario B, Figure \ref{fig:hourlycomparison} plots the total amount of charging per hour for both scenarios. The charging duration on the y-axis is the sum across all three chargers considered, so the maximum possible value would be 180 minutes. We can see that in Scenario A, charging is fairly evenly distributed throughout the day, with decreased charging at the earliest and latest parts of the day when fewer buses are in service. When we introduce delays by increasing trip durations between 7:00 and 9:00 in Scenario B, our method shifts charging activities away from this period---particularly between 8:00 and 11:00, when buses were likely to be running behind schedule due to a previous delay not caused by charging. Accordingly, there is more charging activity in the early morning from 6:00 to 7:00 as well as for most of the afternoon; starting around 8:00 p.m., the charging pattern looks nearly identical to that of Scenario A.

Figure \ref{fig:burientimeline} plots the timelines of charger utilization in Scenarios A and B in the same style as Figure \ref{fig:simpletimeline}. Figure \ref{fig:burientimeline} only shows buses that charge at the Burien Transit Center, which is the busiest of the three chargers. As before, gray blocks on the timeline represent on-time service trips, orange blocks represent delayed trips, and blue blocks indicate charging in Burien. The period from 7:00 to 9:00 a.m., when trips take longer than expected and force delays in Scenario B, is highlighted with a yellow background in the timelines for both scenarios.

Comparing both scenarios shown in Figure \ref{fig:burientimeline}, we can see that many buses that charged in the morning commute hours in Scenario A instead charge earlier or later in the day under Scenario B, in agreement with Figure \ref{fig:hourlycomparison}. In particular, we see for Scenario B that buses generally avoid charging prior to a trip that is delayed, in order to avoid prolonging those delays. The charger is used scarcely between about 8:00 and 11:00 a.m. to allow buses to catch back up to schedule, after which frequent charging resumes. Figure \ref{fig:burientimeline} also shows that most buses charge multiple times throughout the day; in fact, in both scenarios, the median number of times for a bus to charge is 3. Also note that in both scenarios, buses rarely charge for the full amount of layover time available to them. This result suggests in our test cases, individual blocks have plenty of layover time scheduled to allow for recharging, and queuing is the main potential source of delays.

\section{Discussion and Conclusion}
\label{sec:conclusion}
This study developed a novel model for scheduling within-day electric bus fleet recharging as well as exact and heuristic methods for its solution. Our method is unique in that it exactly tracks charger utilization and queuing with a focus on schedule adherence, but with a particularly flexible approach. Our model precisely sets delays and propagates them across trips with the objective of minimizing total delays. Because of this precise yet flexible approach, our model still produces actionable results when operating conditions such as traffic delays or high energy consumption make it impossible to meet charging needs without delaying any trips. This is unique in the literature, as existing methods tend to constrain charging to take place during scheduled layover time and are not naturally extendable to scenarios where insufficient time is available.

Our exact solution method based on Combinatorial Benders decomposition outperforms a state-of-the-art commercial solver on some instances. However, our results show that this approach is still not efficient enough to be suitable for large real-world transit networks under demanding conditions. Although this method can provide a proof of optimality for some instances that provides theoretical value, it has considerable limitations. In particular, our usage of Combinatorial Benders for a problem in which the objective value is set solely by the subproblem is not ideal.

However, our Select--Sequence--Schedule heuristic algorithm that was inspired by the CB approach shows excellent results on our test instances. On our smaller test network instances in Section \ref{sec:simplenetwork}, the 3S method found good feasible solutions orders of magnitude faster than exact methods and its optimality gap was always below 8\%. The two test scenarios on the King County Metro network from Section \ref{sec:metronetwork} demonstrated how the heuristic could quickly find a zero-delay solution under favorable conditions, as well as its ability to respond to exogenous delays by adjusting the recharging plan. Its strong performance in finding good feasible solutions suggests that 3S could also be run repeatedly for the purpose of real-time scheduling, taking into account up-to-the-minute information on delays, traffic conditions, and energy consumption predictions. This approach could help agencies respond to a variety of disruptions such as charger breakdowns while minimizing the ultimate service impacts.

Our method is also straightforward to adapt to support other objectives besides delay minimization. For instance, the objective function \eqref{eq:mipobj2} can be adapted to incorporate other concerns such as time-dependent energy costs rather than focusing purely on delays. Agencies may also wish to prioritize delays on some trips more than others based on expected passenger volumes or equity concerns. Such features can be added to our model and solution methods with minimal changes required. One more practical extension to our work is to relax Assumption \ref{assume:singlecharger} that only one charger is located at each terminal. This change can be accommodated by adding a constraint that only one charger is used per trip and adjusting Phase 1 of the 3S heuristic to ensure this constraint is respected.

\section*{Acknowledgement}
The work of the first author is supported by the National Science Foundation (NSF) Graduate Research Fellowship Program under Grant No. DGE-1762114. This research is also partially supported by a grant from the Pacific Northwest Transportation Consortium (PacTrans) funded by US Department of Transportation (USDOT). Any opinions, findings, and conclusions or recommendations expressed in this paper are those of the authors and do not necessarily reflect the views of the NSF, USDOT, or PacTrans. This work was in part supported by the Transport Area of Advance within Chalmers University of Technology.

\bibliographystyle{elsarticle-harv}
\bibliography{TREDraft}
\end{document}